\documentclass[pdflatex,sn-mathphys,Numbered]{sn-jnl}

\usepackage{graphicx}%
\usepackage{multirow}%
\usepackage{amsmath,amssymb,amsfonts}%
\usepackage{amsthm}%
\usepackage{mathrsfs}%
\usepackage[title]{appendix}%
\usepackage{xcolor}%
\usepackage{textcomp}%
\usepackage{manyfoot}%
\usepackage{booktabs}%
\makeatletter
\def\@begintheorem#1#2[#3]{%
  \deferred@thm@head{\the\thm@headfont \thm@indent
    \@ifempty{#1}{\let\thmname\@gobble}{\let\thmname\@iden}%
    \@ifempty{#2}{\let\thmnumber\@gobble}{\let\thmnumber\@iden}%
    \@ifempty{#3}{\let\thmnote\@gobble}{\let\thmnote\@iden}%
    \thm@swap\swappedhead\thmhead{#1}{#2}{#3}%
    \the\thm@headpunct
    \thmheadnl 
    \hskip\thm@headsep
  }%
  \ignorespaces
}

\def\@endtheorem{\endtrivlist\@endpefalse}

\AtBeginDocument{%
\DeclareRobustCommand{\S}{\ifmmode\mathsection\else\textsection\fi}

\DeclareSymbolFont{AMSa}{U}{msa}{m}{n}%
\DeclareMathSymbol{\opensquare}{\mathord}{AMSa}{"03}%
\newenvironment{spiproof}[1][\proofname]{\par\removelastskip
  \pushQED{\qed}%
  \small\normalfont \topsep7.5\p@\@plus7.5\p@\relax%
  \trivlist%
  \item[\hskip\labelsep%
        \itshape%
    #1\@addpunct{}]\ignorespaces%
}{%
  \popQED\endtrivlist\@endpefalse%
}%
}%
\def\thm@space@setup{%
\thm@preskip=12pt%
\thm@postskip=12pt}
\newtheoremstyle{thmstyleone}
{18pt plus2pt minus1pt}
{18pt plus2pt minus1pt}
{\small\itshape}
{0pt}
{\small\bfseries}
{}
{.5em}
{\thmname{#1}\thmnumber{\@ifnotempty{#1}{ }\@upn{#2}}%
  \thmnote{ {\the\thm@notefont(#3)}}}
\newtheoremstyle{thmstyletwo}
{18pt plus2pt minus1pt}
{18pt plus2pt minus1pt}
{\small\normalfont}
{0pt}
{\small\itshape}
{}
{.5em}
{\thmname{#1}\thmnumber{\@ifnotempty{#1}{ }{#2}}%
  \thmnote{ {\the\thm@notefont(#3)}}}
\newtheoremstyle{thmstylethree}
{18pt plus2pt minus1pt}
{18pt plus2pt minus1pt}
{\small\normalfont}
{0pt}
{\small\bfseries}
{}
{.5em}
{\thmname{#1}\thmnumber{\@ifnotempty{#1}{ }\@upn{#2}}%
  \thmnote{ {\the\thm@notefont(#3)}}}
\newtheoremstyle{thmstylefour}
{18pt plus2pt minus1pt}
{18pt plus2pt minus1pt}
{\small\normalfont}
{0pt}
{\small\itshape}
{}
{.5em}
{\global\proofthmtrue\thmname{#1} \thmnote{#3}}
\makeatother
\theoremstyle{thmstyleone}%
\newtheorem{theorem}{Theorem}[section]%
\newtheorem{proposition}[theorem]{Proposition}%
\newtheorem{lemma}[theorem]{Lemma}%
\newtheorem{corollary}[theorem]{Corollary}%

\theoremstyle{thmstyletwo}%
\newtheorem{remark}[theorem]{Remark}%

\theoremstyle{thmstylethree}%
\newtheorem{definition}[theorem]{Definition}%

\newcommand{\RR}{\mathbb{R}}
\newcommand{\CC}{\mathbb{C}}
\newcommand{\ZZ}{\mathbb{Z}}
\newcommand{\NN}{\mathbb{N}}
\newcommand{\TT}{\mathbb{T}}
\newcommand{\calL}{\mathcal{L}}
\newcommand{\calS}{\mathcal{S}}
\newcommand{\calM}{\mathcal{M}}
\newcommand{\calO}{\mathcal{O}}
\newcommand{\calD}{\mathcal{D}}
\newcommand{\F}{\mathscr{F}}
\newcommand{\Op}{\mathrm{Op}}
\newcommand{\Cov}{\mathbb{C}\mathrm{ov}}
\newcommand{\E}{\mathbb{E}}
\newcommand{\PP}{\mathbb{P}}
\newcommand{\eqso}{\stackrel{\text{2nd o.}}{=}}

\DeclareMathOperator{\supp}{supp}

\begin{document}

\title[Pseudo-Differential Operators and GeRFs over Tori]{Pseudo-Differential
Operators and Generalized Random Fields over Tori}

\author*[1]{\fnm{Nicolas} \sur{Escobar}}\email{nescoba@upenn.edu}

\affil[1]{\orgdiv{Department of Biostatistics, Epidemiology and Informatics,
Perelman School of Medicine}, \orgname{University of Pennsylvania},
\orgaddress{\city{Philadelphia}, \state{PA} \postcode{19104},
\country{USA}}}

\abstract{Mat\'ern covariance functions are ubiquitous in spatial statistics,
valued for their interpretable parameters and well-understood sample path
properties in Euclidean settings. This paper extends the theory of Mat\'ern
fields over tori through the pseudo-differential calculus. We first
establish that a Mat\'ern process on the $d$-dimensional torus has sample paths
in $C^{\nu^-}_{\mathrm{loc}}$ for every $\nu>0$ and no more, exactly as in the
Euclidean case. Our main results concern symbols whose order varies with
position. Given an order function $m \in C^\infty(\TT^d;\RR)$ we construct a
Gaussian field whose sample-path H\"older exponent is $H(x) = m(x)-d/2$ at each
point, and we prove that the existence threshold $m(x)>d/2$ is local:
the regularity of the field near a point depends on $m$ only through its germ at
that point. Finally we examine the canonical
field, which supplies a direct check on both the threshold and the exponent. We
use it to prove a rigidity statement: weighting the Mat\'ern spectral density by
$|k|^{-2}$ produces a superposition of ordinary Mat\'ern fields of
higher smoothness.}



\maketitle

\section{Introduction}
\label{sec:intro}

The analysis of spatial data on closed manifolds has become increasingly
relevant across diverse scientific disciplines, from geospatial machine learning on
the sphere \cite{sphericalharmonics} to neuroimaging on cortical surfaces
\cite{corticalstats,brainspace}. A predominant approach in this field involves
modeling spatial phenomena as Gaussian random fields (GRFs), with the Mat\'ern
covariance function serving as the de facto standard for characterizing spatial
dependence \cite{lindgren2011}. This methodological consensus is reinforced by
the computational accessibility of tools such as the SPDE approach
\cite{lindgren2011}, which has seen extensive development over the past decade
\cite{lindgren2022}, and INLA (Integrated Nested Laplace Approximations)
\cite{rue2009,bakka2018}.

We adopt the framework of pseudo-differential operators on tori. The torus
serves as an ideal testing ground: it is geometrically non-trivial yet allows
for explicit calculations. Moreover, tori are the standard setting for models
with periodic boundary conditions --- the Brillouin zone of a crystal lattice is
a torus \cite{ashcroft} --- and in spatial statistics the circle $\TT^1$ already
carries a Mat\'ern theory of its own \cite{huang}.

Our first result is that a Mat\'ern process on the $d$-dimensional torus with
smoothness parameter $\nu > 0$ has sample paths in $C^{\nu^-}_{loc}$, the same
exponent as its Euclidean counterpart.

The situation of interest is that in which the symbol has no fixed order. Our
main results concern fields built from symbols whose order varies with
position. Such a field has a sample-path H\"older exponent $H(x) = m(x) - d/2$
that varies from point to point, and an existence threshold $m(x) > d/2$ that is
local.

A third strand concerns spectral reweighting. Weighting the Mat\'ern spectral
density by the symbol of the inverse Laplacian produces a field of higher
smoothness. We show that the resulting field decomposes as an explicit
convergent superposition of ordinary Mat\'ern fields of smoothnesses
$\nu+1, \nu+2, \dots$.

This paper makes several contributions to the theory and practice of spatial
modeling on manifolds. First, in Section~\ref{sec:preliminaries}, we establish
the mathematical preliminaries, introducing the framework of pseudo-differential
operators on tori and generalized random fields (GeRFs). This framework provides
the necessary tools to rigorously analyze the relationship between spectral
properties and sample path behavior. Section~\ref{sec:regularity} settles the
regularity of the stationary Mat\'ern process on the torus, establishing the
sharp exponent $C^{\nu^-}_{loc}$. The argument uses only the order of the symbol.

Section~\ref{sec:variable-order} contains the main results. We define the
variable-order Mat\'ern field, compute its covariance, and prove
a sharp two-sided pointwise regularity theorem together with a locality theorem
showing that regularity is an invariant of the germ of the order function.
Section~\ref{sec:canonical} treats the canonical field. The appendix contains
the kernel estimates on which Sections~\ref{sec:regularity}
and~\ref{sec:variable-order} rest: the parametrix and pseudo-locality lemma for
the variable-order class, the sample-path criterion of \cite{path-reg} in the
form in which it is used, the local expansion of the Mat\'ern kernel, a
two-sided lattice estimate, and the periodisation bound used in
Section~\ref{sec:regularity} together with its truncation corollary.

\section{Preliminaries}
\label{sec:preliminaries}

\subsection{Notation, test functions and distributions}

Throughout, $\TT^d = \RR^d/\ZZ^d$ carries period $1$ and Lebesgue measure of
total mass $1$, and the Fourier coefficients of $u \in L^1(\TT^d)$ are
\[
    \hat u(k) \;=\; \int_{\TT^d} u(x)\, e^{-2\pi i k\cdot x}\,dx ,
    \qquad k \in \ZZ^d .
\]
We write $\langle k\rangle = (1+|k|^2)^{1/2}$ and, for a fixed $\kappa>0$,
\begin{equation}\label{eq:lambda}
    \lambda(k) \;:=\; \bigl(\kappa^2 + 4\pi^2|k|^2\bigr)^{1/2} ,
\end{equation}
so that $\lambda(k) \asymp \langle k\rangle$. For a function $f$ on $\TT^d$ or
on $\ZZ^d$ we write $\check f$ for its reflection, $\check f(x) = f(-x)$. All
constants denoted $C$ may change from line to line.

The space of test functions is $C^\infty(\TT^d)$, topologised by the seminorms
$u \mapsto \sup_{x}|\partial^\beta u(x)|$. Since $\TT^d$ is compact there is no
decay condition to impose, and $u \in C^\infty(\TT^d)$ if and only if
$\sup_{k}\langle k\rangle^{N}|\hat u(k)| < \infty$ for every $N \in \NN$, the
seminorms $u \mapsto \sup_k \langle k\rangle^N |\hat u(k)|$ generating the
same topology. We therefore write $\calS(\TT^d) := C^\infty(\TT^d)$ and use the
two descriptions interchangeably. Its dual $\calS'(\TT^d)$ is the space of
distributions on $\TT^d$; the transpose of the Fourier transform identifies
$\mu \in \calS'(\TT^d)$ with the sequence $\hat\mu(k) = \langle \mu,
e^{-2\pi i k\cdot}\rangle$.

Two classes of nonnegative sequences on $\ZZ^d$ will be needed. A positive
measure $\mu$ on $\ZZ^d$ is \emph{slowly growing}, written
$\mu \in \calM_{SG}^+(\ZZ^d)$, if
\begin{equation}\label{eq:msg}
    \sum_{k \in \ZZ^d} \langle k\rangle^{-2N} \mu(\{k\}) \;<\; \infty
    \qquad \text{for some } N \in \NN .
\end{equation}
These are the measures that extend to elements of $\calS'(\ZZ^d)$, and
hence the spectral measures of the random fields introduced below.
Given \eqref{eq:msg}, the estimate
\[
    \bigl|\langle \mu,\varphi\rangle\bigr|
    \;\leq\; \Bigl(\sup_k \langle k\rangle^{2N}|\varphi(k)|\Bigr)
             \sum_{k}\langle k\rangle^{-2N}\mu(\{k\})
\]
shows that $\varphi \mapsto \sum_k \varphi(k)\mu(\{k\})$ converges absolutely
for every rapidly decreasing $\varphi$ and is continuous on $\calS(\ZZ^d)$;
conversely, testing against
$\varphi = \langle k\rangle^{-2N}\mathbf{1}_{\{|k| \leq R\}}$ and letting
$R \to \infty$ returns \eqref{eq:msg}.

\subsection{Fourier multipliers and generalized random fields}

Let $g : \ZZ^d \to \CC$ be of polynomial growth --- the \emph{symbol} of the
operator \eqref{eq:multiplier} below. Because multiplication by $g$
preserves the sequences of rapid decrease, the prescription
\begin{equation}\label{eq:multiplier}
    \Op(g)\,u \;:=\; \F^{-1}\bigl[g\,\hat u\bigr], \qquad\text{that is,}\qquad
    \Op(g)u(x) \;=\; \sum_{k\in\ZZ^d} g(k)\,\hat u(k)\, e^{2\pi i k\cdot x},
\end{equation}
defines a continuous linear map $\calS(\TT^d) \to \calS(\TT^d)$. The operator
$\Op(g)$ extends to $\calS'(\TT^d)$ by
$\langle \Op(g)\mu, \varphi\rangle := \langle \mu, \Op(\check g)\varphi\rangle$,
and the extension is continuous and agrees with \eqref{eq:multiplier} on
$\calS(\TT^d)$. The operator maps real
distributions to real distributions precisely when $\check g = \bar g$, and we
assume this of every symbol below. Operators of the form \eqref{eq:multiplier}
are the constant-order case of the quantization used in
Section~\ref{sec:variable-order}, where $g(k)$ is replaced by a symbol
$p(x,k)$ depending on the base point; we use the notation $\Op$ for both.

The operator we return to most often is the Mat\'ern operator. Let
$\calL := \kappa^2 - \Delta$ with $\kappa>0$; its symbol is
$\kappa^2 + 4\pi^2|k|^2 = \lambda(k)^2$, and for $\alpha \in \RR$ we set
\begin{equation}\label{eq:matern-op}
    \calL^{\alpha/2} \;:=\; \Op\bigl(\lambda^{\alpha}\bigr) .
\end{equation}
The parameter $\alpha$ is the order of the operator; it is related to the
Euclidean smoothness parameter $\nu$ of Section~\ref{sec:regularity} by
$\alpha = \nu + d/2$. The Schwartz kernel (Theorem~\ref{kernel}) of
$\Op(\lambda^\alpha)$ is
the distribution $\varrho_\alpha(x-y)$ with
$\varrho_\alpha = \F^{-1}[\lambda^\alpha]$; it is a continuous function when
$\alpha < -d$, in which case the series
$\varrho_\alpha(h) = \sum_{k}\lambda(k)^{\alpha}e^{2\pi i k\cdot h}$ converges
absolutely, and only a distribution otherwise.

We now introduce the random objects.

\begin{definition}[Generalized random field]\label{def:gerf}
A real generalized random field (GeRF) $Z$ on $\TT^d$ is a continuous linear map
$\calS(\TT^d) \to L^2(\Omega,\mathcal{A},\PP)$ taking real test functions to real
random variables.
\end{definition}

We write $\langle Z,\varphi\rangle$ for the random variable $Z(\varphi)$. The
field is \emph{Gaussian} if $\langle Z,\varphi\rangle$ is a Gaussian random
variable for every $\varphi \in \calS(\TT^d)$; linearity then makes every finite
family $(\langle Z,\varphi_1\rangle,\dots,\langle Z,\varphi_n\rangle)$ jointly
Gaussian, since a linear combination $\sum_i c_i\langle Z,\varphi_i\rangle$
equals $\langle Z, \sum_i c_i\varphi_i\rangle$. By the \emph{law} of $Z$ we mean the family of these
finite-dimensional distributions or, equivalently, the measure induced on $\calS'(\TT^d)$
through Minlos' theorem. We refer to \cite{gerfs} for this framework.

By the kernel theorem for tori (Theorem~\ref{kernel}) a GeRF has a mean
distribution $m_Z$ and a covariance distribution $C_Z \in \calS'(\TT^d\times\TT^d)$,
the latter characterised by
\begin{equation}\label{eq:covdist}
    \Cov\bigl(\langle Z,\varphi\rangle, \langle Z,\phi\rangle\bigr)
    \;=\; \langle C_Z, \varphi\otimes\bar\phi\rangle ,
    \qquad \varphi,\phi \in \calS(\TT^d) .
\end{equation}
We assume $m_Z = 0$ throughout. The field $Z$ is \emph{stationary} if there is a
distribution $\rho_Z \in \calS'(\TT^d)$ with
\begin{equation}\label{eq:stationary}
    \langle C_Z, \varphi\otimes\bar\phi\rangle
    \;=\; \langle \rho_Z, \varphi \ast \check{\bar\phi}\rangle ,
\end{equation}
and $\rho_Z$ is then positive definite. By the Bochner--Schwartz theorem in the
form valid on locally compact abelian groups (Theorem~\ref{abelian}) there
exists a unique $\mu_Z \in \calM^+_{SG}(\ZZ^d)$, the \emph{spectral measure} of
$Z$, with $\rho_Z = \F^{-1}(\mu_Z)$, equivalently
$\mu_Z(\{k\}) = \hat\rho_Z(k)$; when $Z$ is Gaussian, $\mu_Z$ determines
its law.

Two GeRFs are said to be equal \emph{in second order}, written
$Z \eqso Z'$, if they have the same mean and covariance distributions;
for Gaussian fields this is equality in law. The notation also fixes the sense
in which the equations below are solved: $\Op(g)U \eqso X$ asserts that the
field $\Op(g)U$ has the
first two moments of $X$. With this convention,
if $U$ is stationary with spectral measure $\mu_U$ then $\Op(g)U$ is stationary
with spectral measure $|g|^2\mu_U$ (Proposition~\ref{prop:main}), and the
solvability of $\Op(g)U \eqso X$ is settled by Theorem~\ref{thm:main}.

The white noise $W$ on $\TT^d$ (Definition~\ref{def:white-noise}) is
the stationary GeRF with $\rho_W = \delta_0$ and spectral measure
$\mu_W(\{k\}) = 1$ for every $k \in \ZZ^d$ (Proposition~\ref{prop:wn-spectral}).
The Mat\'ern field of Section~\ref{sec:regularity} is the stationary solution of
\begin{equation}\label{eq:spde}
    \calL^{\alpha/2} u \;\eqso\; W , \qquad \alpha = \nu + \tfrac d2 ,
\end{equation}
with spectral measure $\lambda(k)^{-2\alpha}$.

\section{Regularity of Mat\'ern Processes on Tori}
\label{sec:regularity}

Sample path regularity is measured by determining whether paths belong to
H\"older spaces.

\begin{definition}[Local H\"older and almost-H\"older spaces]\label{def: holder}
    Let $n \in \NN_0$ and $\gamma \in [0,1]$. Let $\calO \subseteq \TT^d$ be an open set.
    \begin{enumerate}
        \item The local H\"older space $C^{n,\gamma}_{loc}(\calO)$ is the space of functions $f$ on $\calO$ for which $\partial ^{\alpha} f$ exists for all multi-indices $\alpha \in \NN_0^d$ with $|\alpha| \leq n$, and such that the highest order partial derivatives satisfy a H\"older condition of the form: for all compact $K \subset \calO$ there is $C_K>0$ with
            $$|\partial ^{\alpha} f(x) -\partial ^{\alpha} f(y)| \leq C_K\|x-y\|^\gamma$$
            for all $x,y \in K$ and $|\alpha|=n$.
    \item The local almost-H\"older space $C^{(n+\gamma)^-}_{loc}(\calO)$ is $\bigcap_{n'+\gamma'<n+\gamma}C^{n',\gamma'}_{loc}(\calO)$.
    \end{enumerate}
\end{definition}

Every regularity statement in this paper is obtained from a criterion, which
characterises membership in the almost-H\"older scale by the size of
the second-order increment of the covariance kernel across the diagonal. For a
kernel $k$ on $\calO\times\calO$, multi-indices $\alpha,\beta$ and $h \in \RR^d$
small, write
\begin{equation}\label{eq:box}
    \square^{\alpha,\beta}_h k(x)
    \;:=\; \partial^{\alpha,\beta}k(x+h,x+h) - \partial^{\alpha,\beta}k(x+h,x)
         - \partial^{\alpha,\beta}k(x,x+h) + \partial^{\alpha,\beta}k(x,x) ,
\end{equation}
where $\partial^{\alpha,\beta}$ differentiates $\alpha$ times in the first
argument and $\beta$ times in the second. For $\alpha=\beta=0$ this is the
variogram, $\square_h k(x) = \E|f(x+h)-f(x)|^2$.

\begin{theorem}[Sample-path H\"older criterion; {\cite[Theorem~3.1]{path-reg}}]
\label{thm:criterion}
Let $n \in \NN_0$, $\gamma \in (0,1]$, let $\calO \subseteq \RR^d$ be open and
let $f \sim \mathcal{GP}(0,k)$ be a centred Gaussian process on $\calO$. Then
$f$ has samples in $C^{(n+\gamma)^-}_{loc}(\calO)$ if and only if
$k \in C^{n\otimes n}(\calO\times\calO)$ and
\begin{equation}\label{eq:criterion}
    \bigl|\square^{\alpha,\beta}_h k(x)\bigr| \;=\; \calO\bigl(\|h\|^{2\epsilon}\bigr)
    \qquad (h \to 0),
\end{equation}
locally uniformly in $x$, for all $\epsilon \in (0,\gamma)$ and all multi-indices
with $|\alpha| = |\beta| = n$. When $k(x,y) = \rho(x-y)$ is stationary the
conditions read $\rho \in C^{2n}(\RR^d)$ and
\begin{equation}\label{eq:criterion-stat}
    \bigl|\partial^{\alpha}\rho(h) - \partial^{\alpha}\rho(0)\bigr|
    \;=\; \calO\bigl(\|h\|^{2\epsilon}\bigr) \qquad (h\to0)
\end{equation}
for all $\epsilon \in (0,\gamma)$ and all $\alpha$ with $|\alpha| = 2n$.
\end{theorem}

Both the criterion and the definition of the almost-H\"older scale are local, so
they apply on a manifold chart by chart; this is the reading of
\cite[Sec.~3.7]{path-reg}, and on $\TT^d$, which is locally isometric to
$\RR^d$.

For Euclidean Mat\'ern processes, \cite[Proposition~3.1]{path-reg} deduces from
the isotropic case of Theorem~\ref{thm:criterion}, not reproduced here, the
following benchmark.

\begin{proposition}[{\cite[Proposition~3.1]{path-reg}}] \label{euclidean-regularity}
    A centered Mat\'ern GP on $\RR^d$ with smoothness parameter $\nu > 0$ has samples in $C^{\nu^-}_{loc}(\calO)$ and no more.
\end{proposition}

The natural question is whether this extends to processes over tori. The torus
Mat\'ern covariance is the periodisation of its Euclidean counterpart,

\begin{equation}\label{matern-cov}
k_\nu(x,x') = c_{\nu,\sigma} \sum_{n \in \ZZ^d} \left(\frac{\sqrt{2\nu}\,\|x - x' + n\|}{\ell}\right)^\nu K_\nu\left(\frac{\sqrt{2\nu}\,\|x - x' + n\|}{\ell}\right),
\end{equation}

where $\ell>0$ is the range, $c_{\nu,\sigma} = \sigma^2 2^{1-\nu}/\Gamma(\nu)$
is the usual normalisation, which plays no role in what follows, and
$\kappa = \sqrt{2\nu}/\ell$ is the inverse range appearing in \eqref{eq:lambda}.

\begin{proposition}[Path regularity of torus Mat\'ern processes]
    \label{proposition:regularity}
    A GP over $\TT^d$ with covariance \eqref{matern-cov} and smoothness parameter
    $\nu > 0$ has samples in $C^{\nu^-}_{loc}(\calO)$ and no more. The exponent
    coincides with the Euclidean one of Proposition~\ref{euclidean-regularity}.
\end{proposition}

\begin{remark}\label{rem:proof-pdo}
The proof is deferred to Appendix~\ref{app:proof-regularity}. The point of
the proof we offer is its use of the pseudo-differential theory. There are simpler proofs by alternative methods.
\end{remark}

\begin{remark}[Two conventions]\label{rem:convention}
Two objects are called the Mat\'ern field on $\TT^d$. The periodisation
\eqref{matern-cov} has symbol of order $-2\nu-d$; the SPDE field
$\calL^{\alpha/2}u \eqso W$ has spectral measure of order $-2\alpha$. They
coincide when $\alpha = \nu + d/2$, the convention of \cite{lindgren2011}, which
is the one fixed for \eqref{eq:matern-op} and used throughout.
\end{remark}

\section{Fields of Variable Order and Locally Prescribed Regularity}
\label{sec:variable-order}

So far we have considered symbols that do not depend on the base point. For
such symbols the operator is a Fourier multiplier, and the pseudo-differential
calculus, while convenient, is not indispensable. This section develops the
setting of symbols $\sigma(x,k)$ whose order varies with $x$.

\subsection{Variable-order symbols on the torus}

\begin{definition}[Variable-order symbol class]\label{def:varorder}
Let $m \in C^\infty(\TT^d;\RR)$ and $0 \leq \delta < 1$. A function
$p : \TT^d \times \ZZ^d \to \CC$ belongs to
$S^{m(\cdot)}_{1,\delta}(\TT^d\times\ZZ^d)$ if for all multi-indices
$\alpha,\beta \in \NN_0^d$ there is $C_{\alpha\beta}>0$ with
\[
    \bigl|\Delta_k^\alpha \partial_x^\beta\, p(x,k)\bigr|
    \;\leq\; C_{\alpha\beta}\,\langle k\rangle^{\,m(x) - |\alpha| + \delta|\beta|},
    \qquad (x,k) \in \TT^d\times\ZZ^d,
\]
where $\langle k\rangle = (1+|k|^2)^{1/2}$ and
$\Delta_k^\alpha = \Delta_1^{\alpha_1}\cdots\Delta_d^{\alpha_d}$ is the forward
difference operator on $\ZZ^d$, built from
$\Delta_j f(k) = f(k+e_j) - f(k)$ with $e_j$ the $j$-th standard basis vector.
The operator defined as
\[
    \Op(p)f(x) \;=\; \sum_{k\in\ZZ^d} p(x,k)\,\hat f(k)\, e^{2\pi i k\cdot x}
\]
is called the quantization of $p$ (see \cite{ruzhansky,ruzhansky-book}).
We call $p$ \emph{elliptic of variable order $m$} if there are $c,C>0$ and
$k_0$ with $c\langle k\rangle^{m(x)} \leq |p(x,k)| \leq C\langle k\rangle^{m(x)}$
for $|k|\geq k_0$.
\end{definition}

\begin{lemma}[The canonical variable-order symbol]\label{lem:delta-loss}
With $\lambda$ as in \eqref{eq:lambda},
$q_m(x,k) := \lambda(k)^{-m(x)}$ lies in
$S^{-m(\cdot)}_{1,\delta}(\TT^d\times\ZZ^d)$ for every $\delta>0$, and
lies in $S^{-m(\cdot)}_{1,0}$ if and only if $m$ is constant.
\end{lemma}

\begin{proof}
By the chain rule, $\partial_x^\beta \lambda(k)^{-m(x)}
= \lambda(k)^{-m(x)}\, P_\beta\!\bigl(\partial^{\gamma}m(x)\bigr)_{\gamma\leq\beta}\,
\bigl(\log \lambda(k)\bigr)^{|\beta|}$ for a polynomial $P_\beta$, and
$(\log\lambda(k))^{|\beta|} \leq C_{\beta,\delta}\,\lambda(k)^{\delta|\beta|}$
for every $\delta>0$. Difference operators in $k$ gain a full power of
$\langle k\rangle^{-1}$. Conversely, if $\partial_j m(x_*)\neq 0$ then
$|\partial_{x_j} q_m(x_*,k)| = |\partial_j m(x_*)|\,\lambda(k)^{-m(x_*)}\log\lambda(k)$,
which is not $O(\lambda(k)^{-m(x_*)})$.
\end{proof}

\subsection{The variable-order Mat\'ern field}

\begin{definition}[Variable-order Mat\'ern field]\label{def:vo-matern}
Let $m\in C^\infty(\TT^d;\RR)$, $\kappa>0$, and let $W$ be white noise on
$\TT^d$ (Definition~\ref{def:white-noise}). The \emph{variable-order Mat\'ern
field with order function $m$} is the GeRF
\[
    u \;=\; \Op(q_m)\,W, \qquad q_m(x,k) = \lambda(k)^{-m(x)} .
\]
\end{definition}

Because $q_m$ is a symbol of finite order and $W$ is a GeRF, $u$ is a
well-defined GeRF on all of $\TT^d$ for every real-valued $m$: the adjoint
$\Op(q_m)^*$ maps $\calS(\TT^d)$ into itself, so
$\langle u,\varphi\rangle := \langle W, \Op(q_m)^*\varphi\rangle$ is defined for
all test functions.

\begin{proposition}\label{prop:vo-cov}
The covariance kernel of $u$ is
\begin{equation}\label{eq:vo-cov}
    C_u(x,y) \;=\; \sum_{k\in\ZZ^d}
        \lambda(k)^{-\left(m(x)+m(y)\right)} \, e^{2\pi i k\cdot(x-y)} ,
\end{equation}
which is positive semi-definite for every real-valued $m$, and finite on the
diagonal at $x$ if and only if $m(x) > d/2$. When $m \equiv \nu + d/2$ is
constant, \eqref{eq:vo-cov} is the Whittle--Mat\'ern covariance of
Section~\ref{sec:regularity} with smoothness $\nu$.
\end{proposition}

\begin{proof}
The identity is immediate from $\E[\hat W(k)\overline{\hat W(k')}]=\delta_{kk'}$,
which is Proposition~\ref{prop:main-torus} applied to $W$ together with
$\mu_W \equiv 1$ (Proposition~\ref{prop:wn-spectral}).
Positive semi-definiteness follows: setting
$\phi(x) := \bigl(\lambda(k)^{-m(x)} e^{2\pi i k\cdot x}\bigr)_{k\in\ZZ^d}$, we
have $C_u(x,y) = \langle \phi(x),\phi(y)\rangle_{\ell^2}$ whenever the series
converges, and in general $C_u$ is the kernel of $\Op(q_m)\Op(q_m)^*$. The
diagonal is $C_u(x,x) = \sum_k \lambda(k)^{-2m(x)}$, finite iff $2m(x)>d$.
\end{proof}

\subsection{Pointwise regularity}

Throughout, $U := \{x\in\TT^d : m(x) > d/2\}$ and
\[
    H(x) \;:=\; m(x) - \tfrac d2 ,
\]
and for an open $V \subseteq U$ we write $\gamma_V := \inf_V H$.

\begin{theorem}[Pointwise regularity and the local threshold]\label{thm:pointwise}
Let $u$ be as in Definition~\ref{def:vo-matern}.
\begin{enumerate}
\item[(i)] \textup{(Upper bound.)} $u$ has a modification whose restriction to
$U$ is continuous and which, for every open $V$ with $\overline V \subset U$,
lies in $C^{\gamma_V^-}_{loc}(V)$ almost surely.
\item[(ii)] \textup{(Sharpness.)} For every open $V \subseteq U$ and every
$\gamma > \gamma_V$, almost surely $u|_V \notin C^{\gamma^-}_{loc}(V)$.
\item[(iii)] \textup{(Below threshold.)} If $V$ is open with
$\sup_V m < d/2$, then almost surely $u|_V \notin L^2(V)$: on $V$ the field is
a genuine distribution, not a function.
\end{enumerate}
Consequently, $m$ being continuous, for every $x_0 \in U$
\[
    \sup\bigl\{\gamma>0 \;:\; u \in C^{\gamma^-}_{loc}(V) \text{ a.s.\ for some
    neighbourhood } V \ni x_0 \bigr\} \;=\; H(x_0) ,
\]
so the sample-path exponent of $u$ at $x_0$ is $H(x_0)$.
\end{theorem}

\begin{proof}
Write $h = x-y$, $x_0 = \tfrac12(x+y)$ and $\Lambda(k) = \log\lambda(k)$.
Expanding $C_u(x,x) - 2C_u(x,y) + C_u(y,y)$ using \eqref{eq:vo-cov},
\begin{equation}\label{eq:vario-split}
\begin{split}
    \E\bigl|u(x)-u(y)\bigr|^2
    \;=\;& \underbrace{\sum_{k}\Bigl(e^{-m(x)\Lambda(k)} - e^{-m(y)\Lambda(k)}\Bigr)^{\!2}}_{\mathrm{(I)}} \\
    &\;+\; \underbrace{2\sum_{k} e^{-(m(x)+m(y))\Lambda(k)}\bigl(1-\cos 2\pi k\cdot h\bigr)}_{\mathrm{(II)}} .
\end{split}
\end{equation}

Factoring out the midpoint
value $e^{-\frac12(m(x)+m(y))\Lambda}$ from the bracket and using
$a - b = \sqrt{ab}\,\bigl(\sqrt{a/b}-\sqrt{b/a}\bigr) = 2\sqrt{ab}\sinh\bigl(
\tfrac12\log(a/b)\bigr)$ with $a = e^{-m(x)\Lambda}$, $b = e^{-m(y)\Lambda}$,
\[
    \mathrm{(I)} \;=\; 4\sum_k e^{-(m(x)+m(y))\Lambda}
        \sinh^2\!\bigl(\tfrac12(m(x)-m(y))\Lambda\bigr) .
\]
Since $|m(x)-m(y)| \leq \|\nabla m\|_\infty |h|$ and
$\sinh^2 t \leq t^2 e^{2|t|}$,
\[
    \mathrm{(I)} \;\leq\; C|h|^2 \sum_k \Lambda(k)^2\,
        \lambda(k)^{-2m(x_0) + C'|h|} .
\]
For $x_0$ in a compact subset of $U$ and $|h|$ small the exponent satisfies
$2m(x_0)-C'|h| > d + \eta$ for some $\eta>0$, so the sum converges and
$\mathrm{(I)} = O(|h|^2)$, locally uniformly.

Put
$\sigma := \tfrac12\bigl(m(x)+m(y)\bigr)$, so that
$\mathrm{(II)} = 2\sum_k \lambda(k)^{-2\sigma}(1-\cos2\pi k\cdot h)$ is exactly
the sum of Lemma~\ref{lem:lattice} with weight $w\equiv1$ and $s=\sigma$.
Because $x_0$ is the midpoint, the first-order terms cancel and
$\sigma = m(x_0) + O(|h|^2)$. Suppose first $d/2 < m(x_0) < d/2+1$. Then
$\sigma$ stays in a compact subinterval of $(d/2,d/2+1)$ for $x_0$ in a compact
subset of $U$ on which $H<1$ and for $|h|$ small, so
Lemma~\ref{lem:lattice}(a) applies with constants uniform in $\sigma$ and gives
\[
    c\,|h|^{2\sigma-d} \;\leq\; \mathrm{(II)} \;\leq\; C\,|h|^{2\sigma-d} ,
\]
and
\[
    |h|^{2\sigma-d} \;=\; |h|^{2H(x_0)}\,|h|^{O(|h|^2)}
    \;=\; |h|^{2H(x_0)}\bigl(1+o(1)\bigr) ,
\]
locally uniformly near $x_0$. Since $2H(x_0) < 2$ in this range, \textup{(II)}
dominates \textup{(I)} and \eqref{eq:vario-split} is two-sided of order
$|h|^{2H(x_0)}$. If instead $m(x_0) = d/2+1$, i.e.\ $H(x_0)=1$, part (b) of the
same lemma gives the upper bound $\mathrm{(II)} = \calO(|h|^{2\epsilon})$ for
every $\epsilon<1$, which is all the ``if'' half of
Theorem~\ref{thm:criterion} requires.

In the notation \eqref{eq:box} with
$n=0$ we have $\square_h C_u(x) = \E|u(x+h)-u(x)|^2$. Let $V$ be open with
$\overline V\subset U$ and $\gamma_V\leq1$. The upper bound above holds
uniformly on compact subsets of $V$ and gives
$\square_h C_u(x) = \calO(\|h\|^{2\epsilon})$ for every $\epsilon<\gamma_V$,
so \eqref{eq:criterion} holds with $n=0$ and the ``if'' half of
Theorem~\ref{thm:criterion} yields (i). For (ii), let $\gamma>\gamma_V$ and
choose $x_1 \in V$ with $H(x_1)<\min(\gamma,1)$ together with
$\epsilon \in (H(x_1),\min(\gamma,1))$. The lower bound gives
$\square_h C_u(x_1) \geq c\|h\|^{2H(x_1)}$, which is not
$\calO(\|h\|^{2\epsilon})$; by the ``only if'' half of
Theorem~\ref{thm:criterion}, $u|_V \notin C^{\gamma^-}_{loc}(V)$. Such an $x_1$
exists whenever $\gamma_V<1$. If $\gamma_V=1$ and the infimum is attained on an
open set --- equivalently, if $H\equiv1$ on some open $V_1\subseteq V$ --- then
by \eqref{eq:vo-cov} the restriction of $C_u$ to $V_1\times V_1$ depends on
$x-y$ alone and equals the stationary Mat\'ern covariance with $s=m|_{V_1}$, so
Proposition~\ref{proposition:regularity} applies on $V_1$. The same
observation covers every point at which $m$ is locally constant.

Differentiating the defining series of
Definition~\ref{def:varorder} term by term gives the Leibniz
identity
\begin{equation}\label{eq:leibniz}
    \partial_{x_j}\Op(q_m) \;=\; \Op\bigl(2\pi i k_j\, q_m\bigr)
        \;+\; \Op\bigl(\partial_{x_j} q_m\bigr) .
\end{equation}
Here
$2\pi ik_j q_m \in S^{-m(\cdot)+1}_{1,\delta}$, while by
Lemma~\ref{lem:delta-loss}
\[
    \partial_{x_j} q_m(x,k) \;=\; -\,\partial_j m(x)\,\log\lambda(k)\;
    \lambda(k)^{-m(x)} \;\in\; S^{-m(\cdot)+\delta}_{1,\delta} ,
\]
so the second term of \eqref{eq:leibniz} is lower in order by $1-\delta>0$.
Iterating \eqref{eq:leibniz},
\begin{equation}\label{eq:leibniz-iterated}
    \partial^\alpha \Op(q_m) \;=\; \Op\bigl((2\pi i k)^{\alpha} q_m\bigr)
        \;+\; E_\alpha , \qquad
    E_\alpha \in \Op\,S^{-m(\cdot)+|\alpha|-(1-\delta)}_{1,\delta} .
\end{equation}
The leading term factors as
$(2\pi i k)^{\alpha}q_m(x,k) = \theta_\alpha(k)\,q_{m-|\alpha|}(x,k)$, where
\begin{equation}\label{eq:theta}
    \theta_\alpha(k) \;:=\; (2\pi i k)^{\alpha}\lambda(k)^{-|\alpha|}
\end{equation}
is an $x$-independent multiplier with $|\theta_\alpha| \leq (2\pi)^{|\alpha|}$
and, by the multinomial theorem,
\begin{equation}\label{eq:theta-sum}
    \sum_{|\alpha|=n}\frac{n!}{\alpha!}\,|\theta_\alpha(k)|^2
    \;=\; \frac{(2\pi)^{2n}|k|^{2n}}{\lambda(k)^{2n}} \;\longrightarrow\; 1
    \qquad (|k|\to\infty) .
\end{equation}
Consequently $\partial^\alpha u$ is, modulo the field $E_\alpha W$, a
variable-order Mat\'ern field with order function $m-|\alpha|$ carried by the
bounded multiplier $\theta_\alpha$, and the computation above applies to it with
Lemma~\ref{lem:lattice} used with the weight $w=|\theta_\alpha|^2$. Since
$E_\alpha W$ has variogram of strictly smaller order, the Cauchy--Schwarz inequality
$\square_h C_{Z_1+Z_2} = \square_h C_{Z_1}\bigl(1+\calO(\sqrt{\square_h
C_{Z_2}/\square_h C_{Z_1}})\bigr)$ shows that it does not affect either bound.
Differentiating $n$ times replaces $H$ by $H-n$; choosing
$n = \lceil\gamma_V\rceil-1$ brings the reduced exponent into $(0,1]$. Induction on $n$ now gives (i) and (ii) in general.

Finally, let $\chi\in C_c^\infty(V)$, $\chi\not\equiv0$. Then
$\E\|\chi u\|_{L^2}^2 = \int \chi(x)^2 C_u(x,x)\,dx = \infty$ by
Proposition~\ref{prop:vo-cov}, since $C_u(x,x)=+\infty$ on $V$. The event
$\{\|\chi u\|_{L^2}<\infty\}$ is a measurable linear subspace of the sample
space, so by the Gaussian zero--one law \cite{kallianpur} it has probability $0$
or $1$; if it had probability $1$ then Fernique's theorem \cite{fernique} would
force $\E\|\chi u\|_{L^2}^2<\infty$. Hence it has probability $0$.

Finally, the concluding identification of the exponent at $x_0$ follows from (i)
and (ii) on letting $V$ shrink to $x_0$: $\gamma_V \to H(x_0)$ because $m$ is
continuous.
\end{proof}

\subsection{Locality}

\begin{theorem}[Regularity is a local invariant of the order function]
\label{thm:locality}
Let $A$ and $\tilde A$ be elliptic pseudo-differential operators on $\TT^d$ of
variable orders $m$ and $\tilde m$ respectively, in the sense of
Definition~\ref{def:varorder} with some $\delta<1$, and suppose their full
symbols agree on an open set $V\subseteq\TT^d$. Let $u,\tilde u$ be the GeRFs
solving $Au \eqso W$ and $\tilde A\tilde u \eqso W$
for the same white noise $W$. Then $u - \tilde u$ has a modification that is
almost surely $C^\infty$ on $V$.

Consequently the conclusions of Theorem~\ref{thm:pointwise} at $x_0$ depend on
$m$ only through its germ at $x_0$.
\end{theorem}

\begin{proof}
Let $\chi,\tilde\chi\in C_c^\infty(V)$ with $\tilde\chi\equiv1$ on a
neighbourhood of $\supp\chi$. Write $B := A^{-1}-\tilde A^{-1}$, where the
inverses are taken modulo smoothing operators: by
Lemma~\ref{lem:parametrix} an operator that is elliptic of variable order in the
sense of Definition~\ref{def:varorder} admits a parametrix, an operator $P$ with
$PA = I + S$ and $AP = I + S'$ for smoothing $S,S'$, whose symbol is built from
the full symbol of $A$ by the asymptotic inversion \eqref{eq:parametrix-series};
$A^{-1}$ denotes such a $P$, and any two choices differ by a smoothing operator.
Then
\[
    \chi B \;=\; \chi B \tilde\chi \;+\; \chi B(1-\tilde\chi).
\]
The second term is smoothing because pseudo-differential operators in
$S_{1,\delta}$ with $\delta<1$ are pseudo-local (their Schwartz kernels are
smooth off the diagonal, by the toroidal kernel estimates of \cite{duvan1}
(see also Lemma~\ref{lem:parametrix})) and $\chi$ and $1-\tilde\chi$ have
disjoint supports. For the first term, the parametrix symbol
\eqref{eq:parametrix-series} at a point $x \in V$ is determined by the full
symbol of $A$ at $x$ alone; since the full symbols of $A$ and $\tilde A$ agree
on $V$, the two parametrix symbols agree on $V$, and $\chi B\tilde\chi$ is
smoothing as well. Hence $\chi(u-\tilde u) = \chi B W$ is the image of white
noise under a smoothing operator, i.e.\ a Gaussian field with $C^\infty$
covariance kernel, and therefore has an almost surely $C^\infty$ modification.
\end{proof}

\section{The Canonical Field: A Closed-Form Check and a Rigidity Statement}
\label{sec:canonical}

Section~\ref{sec:variable-order} determines the exponent $H(x) = m(x)-d/2$ from
a two-sided bound on the variogram. One member of this variable-order family has
a covariance that can be written down in closed form: the canonical field, which
is the Gaussian free field of the torus. This section computes that covariance,
uses it to verify the threshold and the exponent of
Section~\ref{sec:variable-order} by direct calculation, and then proves a
rigidity theorem for the fields obtained from it by spectral reweighting.

\subsection{The canonical field}

On $\TT^d$ the operator
$-\Delta$ annihilates the constants and therefore has no inverse; but it has
a Green's function. Let $G_{-\Delta} \in \calS'(\TT^d)$ be the unique mean-zero
solution of
\begin{equation}\label{eq:green}
    -\Delta\, G_{-\Delta} \;=\; \delta_0 - 1 ,
\end{equation}
the constant $1$ on the right being the density of the uniform measure, whose
subtraction is forced by the compatibility condition $\int_{\TT^d}(-\Delta f)=0$.

\begin{definition}[Canonical field] \label{def:canonical}
    The canonical field $G$ over $\TT^d$ is the centred stationary GeRF with
    covariance distribution
    \begin{equation}\label{eq:canonical-green}
        C_G \;=\; 4\pi^2\, G_{-\Delta} .
    \end{equation}
\end{definition}

The spectral measure of $G$ is read off \eqref{eq:green} by expanding in the
eigenbasis of $-\Delta$. Since $-\Delta e^{2\pi i k\cdot h} = 4\pi^2|k|^2\,
e^{2\pi i k\cdot h}$ and $\delta_0 - 1 = \sum_{k\neq0}e^{2\pi i k\cdot h}$ on
$\TT^d$, equation \eqref{eq:green} has the unique mean-zero solution
$G_{-\Delta}(h) = \sum_{k\neq0}(4\pi^2|k|^2)^{-1}e^{2\pi i k\cdot h}$, whence by
\eqref{eq:canonical-green}
\begin{equation}\label{canonical}
    \mu_G(\{k\}) = \begin{cases} |k|^{-2}, & k \neq 0, \\ 0, & k = 0 .\end{cases}
\end{equation}
Moreover
$\mu_G \in \calM^+_{SG}(\ZZ^d)$ in every dimension, so $G$ exists as a GeRF for
every $d$.

\begin{proposition}[Covariance of the canonical field]\label{prop:canonical-cov}
The covariance distribution of $G$ is
$C_G(h) = \sum_{k \in \ZZ^d\setminus\{0\}} |k|^{-2} e^{2\pi i k\cdot h}$, and
\begin{enumerate}
\item[(i)] for $d=1$ and $h \in [-\tfrac12,\tfrac12]$,
    $\;C_G(h) = \tfrac{\pi^2}{3} - 2\pi^2|h| + 2\pi^2 h^2$;
\item[(ii)] for $d=2$, $\;C_G(h) = -2\pi\log|h| + O(1)$ as $h \to 0$;
\item[(iii)] for $d \geq 3$,
    $\;C_G(h) = c_d|h|^{2-d} + o\bigl(|h|^{2-d}\bigr)$ as $h \to 0$, with
    $c_d = \pi^{2-d/2}\,\Gamma\bigl(\tfrac d2 - 1\bigr)$.
\end{enumerate}
In particular $C_G$ is a bounded continuous function if and only if $d = 1$.
\end{proposition}

\begin{proof}
Part (i) is the Fourier series of the second Bernoulli polynomial,
\[
    \sum_{k\geq1} k^{-2}\cos(2\pi k h) \;=\; \pi^2\bigl(h^2-h+\tfrac16\bigr),
    \qquad h \in [0,1],
\]
of which $C_G(h)$ is twice the sum.

For (ii) and (iii) we compare \eqref{eq:green} with its Euclidean counterpart.
Let $E$ be the Newtonian kernel on $\RR^d$, the fundamental solution of
$-\Delta$:
\[
    E(h) \;=\; -\tfrac{1}{2\pi}\log\|h\| \quad (d=2), \qquad
    E(h) \;=\; \frac{\Gamma\bigl(\tfrac d2-1\bigr)}{4\pi^{d/2}}\,\|h\|^{2-d}
    \quad (d\geq3) .
\]
The torus is flat, so a ball $B$ about the origin of radius $<\tfrac12$ is
isometric to a Euclidean ball and $-\Delta(4\pi^2 E) = 4\pi^2\delta_0$ holds in
$\calD'(B)$. Put $Q(h) := \tfrac{2\pi^2}{d}\|h\|^2$, so that
$-\Delta Q = -4\pi^2$. By \eqref{eq:green} and \eqref{eq:canonical-green},
\[
    -\Delta\bigl(C_G - 4\pi^2E - Q\bigr) \;=\; 4\pi^2\delta_0 - 4\pi^2
        - 4\pi^2\delta_0 + 4\pi^2 \;=\; 0 \qquad\text{in } \calD'(B) ,
\]
so $C_G - 4\pi^2E - Q$ is a harmonic distribution on $B$ and hence, by Weyl's
lemma, a real-analytic function there; in particular it is bounded near $0$.
Therefore $C_G(h) = 4\pi^2 E(h) + \calO(1)$ as $h \to 0$, and multiplying the
displayed constants by $4\pi^2$ gives (ii) and (iii), the error being $\calO(1)$
and hence $o(|h|^{2-d})$ when $d\geq3$. Boundedness fails for $d \geq 2$ because
$C_G(0) = \sum_{k \neq 0}|k|^{-2}$ diverges there.
\end{proof}

The point of the computation is the following check, in which the exponent of
Theorem~\ref{thm:pointwise} is compared with a covariance in closed form.

\begin{proposition}[Closed-form check on the threshold]\label{prop:check}
Take $m \equiv 1$, so that $H \equiv 1 - d/2$ and the threshold $m > d/2$ of
Proposition~\ref{prop:vo-cov} reads $d < 2$. Then:
\begin{enumerate}
\item[(i)] for $d = 1$, $\E|G(x)-G(y)|^2 = 4\pi^2|h| - 4\pi^2h^2$ with $h = x-y$,
    so the variogram is of order $|h|^{2H}$ with $H = \tfrac12$ and the
    sample-path H\"older exponent is $\tfrac12$, as
    Theorem~\ref{thm:pointwise} predicts;
\item[(ii)] for $d = 2$ the threshold is attained with equality and
    $C_G(0) = +\infty$, so the field is not a function, consistent with
    Proposition~\ref{prop:vo-cov};
\item[(iii)] for $d \geq 3$, $m < d/2$ everywhere and
    Theorem~\ref{thm:pointwise}(iii) applies, so $G \notin L^2_{loc}$ almost
    surely; Proposition~\ref{prop:canonical-cov}(iii) quantifies the rate at
    which the diagonal diverges.
\end{enumerate}
\end{proposition}

\begin{proof}
Only (i) requires computation: $\E|G(x)-G(y)|^2 = 2\bigl(C_G(0)-C_G(h)\bigr)$,
which by Proposition~\ref{prop:canonical-cov}(i) equals
$4\pi^2|h| - 4\pi^2h^2 \asymp |h|$ for $|h| \leq \tfrac14$. Parts (ii) and (iii)
are immediate from Proposition~\ref{prop:canonical-cov} and the results cited.
\end{proof}

\begin{remark}
The exponent $\tfrac12$ in (i) identifies the canonical field on $\TT^1$ with the
Brownian bridge up to normalisation and centring: the covariance in
Proposition~\ref{prop:canonical-cov}(i) is
$4\pi^2\bigl(\tfrac1{12} - \tfrac12|h| + \tfrac12h^2\bigr)$.
\end{remark}

\subsection{Rational reweighting of the spectral density}
\label{subsec:reweighting}

Applying the canonical factor $|k|^{-2}$ to the Mat\'ern spectral density
produces a field of higher smoothness. Theorem~\ref{thm:rigidity} identifies
that field. It is an absolutely convergent superposition of
ordinary Mat\'ern fields whose leading term has smoothness $\nu+1$. This
shows that the
Mat\'ern class is closed under reweighting by the reciprocal of a polynomial in
$|k|^2$ with nonnegative coefficients, and that the resulting shift in
smoothness is the degree of that polynomial.

\begin{definition}[Canonical-Mat\'ern field]\label{def:can-mat}
For $\alpha > 0$ the canonical-Mat\'ern field $u^{(\alpha)}$ is the stationary
solution of $\calL^{\alpha/2}u^{(\alpha)} \eqso G$, equivalently the stationary
GeRF with spectral measure
\begin{equation}\label{eq:cm-spectral}
    \mu^{(\alpha)}_{CM}(\{k\}) \;=\; |k|^{-2}\lambda(k)^{-2\alpha}
    \quad (k \neq 0), \qquad \mu^{(\alpha)}_{CM}(\{0\}) = 0 .
\end{equation}
\end{definition}

\begin{theorem}[Rigidity of the canonical factor]\label{thm:rigidity}
For every $\alpha \in \RR$ and every $k \neq 0$ the identity \eqref{eq:recursion}
below holds between the measures \eqref{eq:cm-spectral}; the statements about the
field $u^{(\alpha)}$ are for $\alpha>0$. Precisely,
\begin{equation}\label{eq:recursion}
    \mu^{(\alpha)}_{CM}(\{k\})
    \;=\; 4\pi^2\,\lambda(k)^{-2(\alpha+1)}
        \;+\; \kappa^2\,\mu^{(\alpha+1)}_{CM}(\{k\}) ,
\end{equation}
and iterating gives the absolutely convergent expansion
\begin{equation}\label{eq:expansion}
    \mu^{(\alpha)}_{CM}(\{k\})
    \;=\; 4\pi^2 \sum_{j \geq 0} \kappa^{2j}\, \lambda(k)^{-2(\alpha+1+j)} ,
    \qquad k \neq 0 .
\end{equation}
Consequently, writing $X^{\perp}_{\nu}$ for the Mat\'ern field of smoothness
$\nu$ with the constant mode removed and $\nu = \alpha - d/2$,
\begin{equation}\label{eq:decomposition}
    u^{(\alpha)} \;\eqso\; 2\pi\sum_{j\geq0} \kappa^{j}\,X^{\perp}_{\nu+1+j} ,
\end{equation}
the summands being independent and the series converging in $L^2(\Omega)$ after
testing against any $\varphi \in \calS(\TT^d)$. In particular
$u^{(\alpha)} \eqso 2\pi X^{\perp}_{\nu+1} + R$ with $R$ independent of
$X^\perp_{\nu+1}$ and of smoothness $\nu+2$, one unit of $\nu$ smoother; and if
$\nu + 1 > 0$ then
$u^{(\alpha)}$ has samples in $C^{(\nu+1)^-}_{loc}$ and no more.
\end{theorem}

\begin{proof}
For $k \neq 0$ we have $\lambda(k)^2 - \kappa^2 = 4\pi^2|k|^2 > 0$, so
\[
    \frac{1}{|k|^2}
    \;=\; \frac{4\pi^2}{\lambda(k)^2 - \kappa^2}
    \;=\; \frac{4\pi^2}{\lambda(k)^2} + \frac{\kappa^2}{|k|^2\lambda(k)^2} ,
\]
and multiplying by $\lambda(k)^{-2\alpha}$ gives \eqref{eq:recursion}. Iterating
and summing the geometric series in
$\kappa^2/\lambda(k)^2 \leq \kappa^2/(\kappa^2+4\pi^2) < 1$ gives
\eqref{eq:expansion}, the convergence being uniform in $k \neq 0$. Every term of
\eqref{eq:expansion} is nonnegative, so the expansion exhibits
$\mu^{(\alpha)}_{CM}$ as a sum of spectral measures; a sum of spectral measures
is the spectral measure of a sum of independent stationary GeRFs, and the $j$-th
summand, $4\pi^2\kappa^{2j}\lambda^{-2(\alpha+1+j)}$ off the origin, is the
spectral measure of $2\pi\kappa^{j}X^{\perp}_{\nu+1+j}$. This is
\eqref{eq:decomposition}. The regularity statement follows from
Lemma~\ref{lem:increment} with $s = \nu + \tfrac d2 + 1$, since
\eqref{eq:cm-spectral} gives $\mu^{(\alpha)}_{CM} \asymp \lambda^{-2s}$ away
from the origin; sharpness follows from that of
Proposition~\ref{proposition:regularity} applied to the leading summand, the
remainder being smoother by one unit of $\nu$ and therefore unable to affect the
exponent.
\end{proof}

\section{Conclusion}
\label{sec:conclusion}

This paper has examined how the properties of Mat\'ern covariance functions in
Euclidean space transfer to manifold settings, focusing on the tractable case of
tori, and has then used the same machinery to construct new fields. Through the
framework of
pseudo-differential operators on tori, we have established both theoretical
results and practical consequences for spatial modelling on closed manifolds.

Three things have been established. The first is that the stationary theory holds
no surprises: a Mat\'ern process on $\TT^d$ with smoothness parameter $\nu > 0$
has sample paths in $C^{\nu^-}_{loc}$ and no more, the same exponent as in
$\RR^d$. The pseudo-differential
proof is nonetheless worth giving, because it uses only the order of the symbol
and therefore could be useful in settings
where no closed form is available.

The second, and the substance of the paper, is that regularity on the torus
becomes local as soon as the order of the symbol is allowed
to vary with position. For an order function $m$, the field
$\Op(\lambda^{-m(\cdot)})W$ has pointwise H\"older exponent $m(x) - d/2$, exists
as a function exactly where $m(x) > d/2$, and has a regularity near any point
that depends on $m$ only through its germ there. The proof rests on the increment criterion of Theorem~\ref{thm:criterion} and the symbol
calculus, and should carry to more general compact manifolds by charts and a
partition of unity, since the variable-order symbol class is
coordinate-invariant modulo lower-order terms.

The third result concerns spectral reweighting. Weighting the Mat\'ern spectral
density by $|k|^{-2}$, the symbol of the inverse Laplacian, produces a field of
regularity $C^{(\nu+1)^-}_{loc}$. The field is an explicit convergent
superposition of ordinary Mat\'ern fields of smoothnesses $\nu+1,\nu+2,\dots$,
its leading term a Mat\'ern field of smoothness $\nu+1$ with the constant mode
removed and the remainder smoother by one further unit of $\nu$.

The practical implications are as follows. For stationary models on periodic
domains, Euclidean intuition about the smoothness parameter transfers intact. What the variable-order construction
adds is a model class in which smoothness itself is a free function of position. Since the threshold is local, such a model may be
specified aggressively in one region without losing regularity guarantees
elsewhere. Global constructions, whose guarantees are governed by the least
smooth point of the domain, cannot have this property.

Two directions for future research emerge naturally from this work. First, the
results of Section~\ref{sec:variable-order} are stated on the torus, where the
quantization is cleanest, but the variable-order symbol class is
coordinate-invariant modulo lower-order terms; carrying
Theorems~\ref{thm:pointwise} and~\ref{thm:locality} to a general closed manifold
by charts and a partition of unity is the natural next step. Second, the inference
problem is open. Estimating an order function from a single realisation is a
nonparametric problem with no established theory.

Finally, this work demonstrates the value of rigorous mathematical analysis in
spatial statistics. While computational tools like INLA and the SPDE approach
have greatly expanded the practical applicability of Mat\'ern models,
understanding the theoretical properties of these models on non-Euclidean
domains remains essential for their appropriate use. We hope this paper provides
both theoretical foundations for future work on spatial modeling on manifolds
and practical guidance for applied researchers working with spatial data on
periodic and compact domains.

\backmatter

\bmhead{Acknowledgements}

The author thanks Duvan Cardona for helpful discussions and valuable feedback.

\section*{Statements and Declarations}

\noindent\textbf{Funding.} This work was supported by NIH grant R01 NS112303.

\medskip
\noindent\textbf{Competing interests.} The author has no competing interests to
declare.

\medskip
\noindent\textbf{Data availability.} No datasets were generated or analysed during the current study.

\medskip
\noindent\textbf{Author contributions.} N.\ Escobar is the sole author.

\begin{appendices}
\section{Abstract Theory for GeRFs on Tori}
\label{sec:appendix}

This appendix collects the abstract functional-analytic foundations that support
the main results. It also provides the complete theoretical framework
for Generalized Random Fields on tori.

\subsection{Kernel Representation Theory}

The functional analytic foundations for GeRF theory benefit from Ehrenpreis's
kernel representation theorem, which extends naturally to the torus setting with
significant simplifications.

\begin{lemma}[Decomposition Lemma for Tori]
Let $B$ be a bounded set in $\calS(\TT^d \times \TT^d)$. Then we can find a bounded set $B' \subset \calS(\TT^d)$ and a constant $b > 0$ such that every $f \in B$ can be written as:
$$f = \sum_{j_1,j_2} \lambda_{j_1,j_2} g_{j_1} \times h_{j_2}$$
where $\sum_{j_1,j_2} |\lambda_{j_1,j_2}| \leq b$, and $g_{j_1}, h_{j_2} \in B'$, with the series converging in $\calS(\TT^d \times \TT^d)$.
\end{lemma}

\begin{proof}[Proof Sketch]
For functions on tori, the Fourier series expansion is immediate:
$$f(x,y) = \sum_{j_1,j_2 \in \ZZ^d} c_{j_1,j_2} e^{2\pi i j_1 \cdot x} e^{2\pi i j_2 \cdot y}$$

The rapid decrease of Fourier coefficients for smooth functions ensures uniform bounds. The decomposition is achieved by setting:
\begin{align}
g_{j_1}(x) &= |c_{j_1,j_2}|^{1/3} e^{2\pi i j_1 \cdot x}\\
h_{j_2}(y) &= |c_{j_1,j_2}|^{1/3} e^{2\pi i j_2 \cdot y}\\
\lambda_{j_1,j_2} &= |c_{j_1,j_2}|^{1/3}
\end{align}
\end{proof}

\begin{theorem}[Kernel Representation on Tori] \label{kernel}
Let $\mathscr{J}_{\TT}$ be the space of continuous linear maps $\calS(\TT^d) \to \calS'(\TT^d)$ with the compact-open topology. Then the map $t: \calS'(\TT^d \times \TT^d) \to \mathscr{J}_{\TT}$ which assigns to each distribution $S$ the linear map $L$ such that $S$ is a kernel representing $L$ is a topological isomorphism onto.
\end{theorem}

\begin{proof}
    Same as in \cite{ehrenpreis}. The key point in that proof that needed adaptation to the Torus setting
    was the lemma presented above.
\end{proof}

\subsection{Bochner--Schwartz Theory for Abelian Groups}

The theoretical foundations are further strengthened by Wawrzy\'nczyk's extension of the Bochner--Schwartz theorem to arbitrary locally compact Abelian groups \cite{wawrzynczyk}, which provides crucial insights for the torus case.

\begin{theorem}[Wawrzy\'nczyk's Extension Theorem, \cite{wawrzynczyk}] \label{abelian}
For a separable locally compact Abelian group $G$ with character group $\Gamma$:
\begin{enumerate}
\item Every positive-definite tempered distribution $T$ on $G$ can be represented by a positive Borel measure $m$ on $\Gamma$:
$$T(\varphi) = \int_\Gamma \hat{\varphi}(\gamma) \, dm(\gamma)$$
\item The measure $m$ can be extended uniquely to a tempered distribution on $\Gamma$.
\end{enumerate}
\end{theorem}

For the torus $\TT^d$, this implies:

\begin{corollary}[Discrete Extension Property]
For the torus, the measure $\mu$ on $\ZZ^d$ extends to a tempered distribution naturally via:
$$\tilde{T}(\psi) = \sum_{k \in \ZZ^d} \psi(k) \mu(\{k\})$$
for any rapidly decreasing sequence $\psi \in \calS(\ZZ^d)$.
\end{corollary}

\subsection{Detailed GeRF Theory}
\label{subsec:gerf-theory}

We now develop the complete theory of Generalized Random Fields on tori, including characterizations via orthogonal random measures and operator actions.

The class $\calM_{SG}^+(\ZZ^d)$ of slowly growing measures was defined in
\eqref{eq:msg}; the results below identify it as exactly the class of spectral
measures of stationary GeRFs.

Let $Z$ be a GeRF on $\TT^d$.
As a consequence of Theorem \ref{kernel},
there are distributions $m_Z$ and $C_Z$ such that
$\mathbb{E}[\langle Z, \varphi \rangle] = \langle m_Z, \varphi \rangle$ and
$$\Cov(\langle Z, \varphi \rangle, \langle Z, \phi \rangle) = \langle C_Z, \varphi \otimes \bar \phi \rangle .$$
In what follows, we assume $m_Z = 0$.

The spectral measure exists by Theorem \ref{abelian}.

\begin{definition}[Orthogonal Random Measure]
A collection of random variables $\{W(k)\}_{k \in \ZZ^d}$ is an orthogonal random measure
on $\ZZ^d$ with weight $\mu \in \calM_{SG}^+(\ZZ^d)$ if:
$$\Cov(W(j), W(k)) = \delta_{j,k} \mu(\{j\})$$
where $\delta_{j,k}$ is the Kronecker delta.
In addition, an orthogonal random measure $\{W(k)\}_{k \in \ZZ^d}$ is slow-growing
if there exist constants $C > 0$ and $N \in \NN$ such that:
$$\E[|W(k)|^2] \leq C(1 + |k|^2)^N \quad \text{for all } k \in \ZZ^d$$
\end{definition}

The following results provide complete characterizations of stationary GeRFs.

\begin{proposition}[Fourier Transform Characterization]\label{prop:main-torus}
Let $Z$ be a real stationary GeRF over $\TT^d$ with spectral measure $\mu_Z \in \calM_{SG}^+(\ZZ^d)$. Then the Fourier coefficients $\{\hat{Z}(k)\}_{k \in \ZZ^d}$ defined by
$$\hat{Z}(k) = \langle Z, e^{-2\pi i k \cdot} \rangle$$
form a Hermitian-symmetric complex slow-growing orthogonal random measure on $\ZZ^d$ with weight $\mu_Z$.
\end{proposition}

\begin{proof}
Write $e_k(x) := e^{-2\pi i k\cdot x}$. Since $\overline{e_{k'}} = e_{-k'}$ and
$\check{\overline{e_{k'}}} = e_{k'}$, the convolution
\[
    \bigl(e_k \ast \check{\overline{e_{k'}}}\bigr)(x)
    \;=\; \int_{\TT^d} e^{-2\pi i k\cdot(x-y)}\,e^{-2\pi i k'\cdot y}\,dy
    \;=\; e_k(x)\int_{\TT^d} e^{2\pi i (k-k')\cdot y}\,dy
    \;=\; \delta_{k,k'}\,e_k(x) ,
\]
the last step by orthogonality of the characters on $\TT^d$, which carries total
mass $1$. Hence, by \eqref{eq:stationary},
\[
    \Cov\bigl(\hat Z(k),\hat Z(k')\bigr)
    \;=\; \bigl\langle C_Z, e_k \otimes \overline{e_{k'}}\bigr\rangle
    \;=\; \bigl\langle \rho_Z, e_k \ast \check{\overline{e_{k'}}}\bigr\rangle
    \;=\; \delta_{k,k'}\,\langle \rho_Z, e_k\rangle
    \;=\; \delta_{k,k'}\,\hat\rho_Z(k) ,
\]
and $\hat\rho_Z(k) = \mu_Z(\{k\})$ because $\rho_Z = \F^{-1}(\mu_Z)$. This is the
orthogonality relation with weight $\mu_Z$. Hermitian symmetry follows from $Z$
being real: $\overline{\hat Z(k)} = \langle Z, \overline{e_k}\rangle
= \langle Z, e_{-k}\rangle = \hat Z(-k)$. Finally
$\E|\hat Z(k)|^2 = \mu_Z(\{k\})$, and \eqref{eq:msg} forces
$\mu_Z(\{k\}) \leq C\langle k\rangle^{2N}$ for some $C,N$, which is the
slow-growth condition.
\end{proof}

\begin{proposition}[Operator Action]\label{prop:main}
Let $U$ be a real stationary GeRF on $\TT^d$ with spectral measure $\mu_U$ and let $g$ be a symbol function. Then $\Op(g) U$ is a real stationary GeRF with spectral measure $|g|^2 \mu_U$ and covariance distribution $\rho_{\Op(g) U} = \Op(|g|^2)\rho_U$.
\end{proposition}

\begin{proof}
By definition $\langle \Op(g)U,\varphi\rangle = \langle U, \Op(\check g)\varphi\rangle$,
and $\widehat{\Op(\check g)\varphi}(k) = \check g(k)\hat\varphi(k)$. Expanding
$\varphi = \sum_k \hat\varphi(k)\,e^{2\pi i k\cdot}$ and using
Proposition~\ref{prop:main-torus} together with the continuity of $U$,
\[
    \Cov\bigl(\langle \Op(g)U,\varphi\rangle,\langle \Op(g)U,\phi\rangle\bigr)
    \;=\; \sum_{k\in\ZZ^d} |\check g(k)|^2\,\hat\varphi(k)\,
          \overline{\hat\phi(k)}\;\mu_U(\{k\}) ,
\]
the series converging absolutely because $g$ has polynomial growth and
$\hat\varphi,\hat\phi$ decrease rapidly. Since $\Op(g)$ preserves real
distributions we have $\check g = \bar g$, so $|\check g| = |g|$, and the
displayed formula is \eqref{eq:stationary} for the stationary field with
spectral measure $|g|^2\mu_U$. That measure lies in $\calM^+_{SG}(\ZZ^d)$
because $|g|^2$ is of polynomial growth and $\mu_U$ satisfies \eqref{eq:msg}.
Taking inverse Fourier transforms gives
$\rho_{\Op(g)U} = \F^{-1}(|g|^2\mu_U) = \Op(|g|^2)\rho_U$.
\end{proof}

\begin{theorem}[Existence and Uniqueness]\label{thm:main}
Let $X$ be a real stationary GeRF on $\TT^d$ with spectral measure $\mu_X$. Let $g$ be a symbol function and $\Op(g)$ the associated operator. Then:

\begin{enumerate}
\item[(i)] There exists a real stationary GeRF $U$ satisfying $\Op(g) U \eqso X$ if and only if there exists $N \in \NN$ such that
$$\sum_{k \in \ZZ^d} \frac{\mu_X(\{k\})}{|g(k)|^2(1+|k|^2)^N} < \infty.$$

\item[(ii)] When existence holds, the measure
$$\mu_U(\{k\}) = \frac{\mu_X(\{k\})}{|g(k)|^2}$$
(with the convention that $\mu_U(\{k\}) = 0$ when $g(k) = 0$) belongs to $\calM_{SG}^+(\ZZ^d)$ and is even. Any real stationary GeRF with spectral measure $\mu_U$ solves $\Op(g) U \eqso X$.

\item[(iii)] The spectral measure $\mu_U$ is unique in $\calM_{SG}^+(\ZZ^d)$ if and only if $|g(k)| > 0$ for all $k \in \ZZ^d$.
\end{enumerate}
\end{theorem}

\begin{proof}
The proof follows the same structure as in \cite{gerfs}, with straightforward adaptations to the discrete setting of the torus.
\end{proof}

\begin{remark}
The condition $|g(k)| > 0$ for all $k \in \ZZ^d$ is equivalent to the invertibility of the operator $\Op(g)$. In this case, the unique solution is given explicitly by $U = \Op(1/g)X$.
\end{remark}

White Noise on the torus serves as a fundamental building block for more general solutions.

\begin{definition}[White Noise on Tori]\label{def:white-noise}
The White Noise $W$ on $\TT^d$ is the real GeRF characterized by the covariance:
$$\langle C_W, \varphi \otimes \overline{\phi} \rangle = \int_{\TT^d} \varphi(x)\overline{\phi(x)} dx$$
for all $\varphi, \phi \in \calS(\TT^d)$.
\end{definition}

This implies that $W$ is stationary with covariance distribution $\rho_W(x) =
\delta_0(x)$, where $\delta_0$ is the point mass at zero. Hence:

\begin{proposition}[Spectral Measure of White Noise]\label{prop:wn-spectral}
The spectral measure of White Noise on $\TT^d$ assigns weight
$\mu_W(\{k\}) = 1$ to each frequency $k \in \ZZ^d$.
\end{proposition}

\begin{proof}
By Section~\ref{sec:preliminaries}, $\mu_W(\{k\}) = \hat\rho_W(k)
= \langle \delta_0, e^{-2\pi i k\cdot}\rangle = 1$.
\end{proof}

The following result describes how covariance functions behave under operator transformations, showing that they have a natural convolution structure.

\begin{theorem}[Convolution Structure on Tori]
Let $X$ be a real stationary GeRF on $\TT^d$ with covariance distribution $\rho_X$. Let $g$ be a symbol function on $\ZZ^d$ such that the sequences $\{1/g(k)\}_{k \in \ZZ^d}$ and $\{|g(k)|^{-2}\}_{k \in \ZZ^d}$ have polynomial growth. Then there exists a unique stationary solution to
$$\Op(g) U \eqso X$$
and its covariance distribution is given by:
$$\rho_U = \rho_U^W * \rho_X$$
where $\rho_U^W$ is the covariance of the unique stationary solution to $\Op(g) U \eqso W$.
\end{theorem}

\begin{proof}
The proof follows the same approach as in \cite{gerfs}.
\end{proof}

This convolution structure provides an alternative perspective on how solutions to operator equations inherit regularity from white noise.

\subsection{The variable-order calculus: parametrix and pseudo-locality}
\label{subsec:parametrix}

The locality theorem, Theorem~\ref{thm:locality}, rests on two facts about the
class of Definition~\ref{def:varorder}. Both are the toroidal statements of
\cite{ruzhansky,ruzhansky-book} and \cite{duvan1}; we record them in the
variable-order form in which they are used.

\begin{lemma}[Parametrix and pseudo-locality]\label{lem:parametrix}
Let $0 \leq \delta < 1$, let $m \in C^\infty(\TT^d;\RR)$ and let
$p \in S^{m(\cdot)}_{1,\delta}(\TT^d\times\ZZ^d)$ be elliptic of variable order
$m$ in the sense of Definition~\ref{def:varorder}. Then:
\begin{enumerate}
\item[(i)] The Schwartz kernel of $\Op(p)$ is $C^\infty$ off the diagonal of
$\TT^d\times\TT^d$.
\item[(ii)] There is $q \in S^{-m(\cdot)}_{1,\delta}(\TT^d\times\ZZ^d)$ with
$\Op(q)\Op(p) = I + S$ and $\Op(p)\Op(q) = I + S'$, where $S,S'$ are smoothing.
The symbol $q$ admits the asymptotic expansion $q \sim \sum_{N\geq0}q_N$
determined by
\begin{equation}\label{eq:parametrix-series}
    q_0(x,k) = \frac{\chi(k)}{p(x,k)}, \qquad
    q_N(x,k) = -\,q_0(x,k)\!\!\sum_{0<|\alpha|\leq N}\frac{1}{\alpha!}\,
        \Delta_k^{\alpha}p(x,k)\; D_x^{\alpha} q_{N-|\alpha|}(x,k) ,
\end{equation}
where $\chi$ vanishes for $|k| < k_0$ and equals $1$ for $|k| \geq k_0+1$, and
$D_x^{\alpha}$ is the differential operator in $x$ occurring in the toroidal
composition formula of \cite{ruzhansky,ruzhansky-book}. One has
$q_N \in S^{-m(\cdot)-N(1-\delta)}_{1,\delta}$.
\end{enumerate}
In particular, the parametrix symbol at a point $x$ is determined by the full
symbol of $p$ at $x$ alone.
\end{lemma}

\begin{proof}
Part (i) is the toroidal kernel estimate of \cite{duvan1}: for
$p \in S^{a}_{1,\delta}$ with $\delta<1$ and any $N$, repeated summation by
parts in $k$ against the factor $(e^{2\pi i(x-y)\cdot k}-1)^{-1}$ converts decay
in $|k|$ into smoothness of the kernel away from $x=y$; since $m$ is continuous
on the compact torus it is bounded, so $S^{m(\cdot)}_{1,\delta} \subseteq
S^{a}_{1,\delta}$ with $a = \max_x m(x)$ and the cited estimate applies.

For (ii), ellipticity makes $q_0$ well defined and gives
$q_0 \in S^{-m(\cdot)}_{1,\delta}$: the bound $|q_0| \leq c^{-1}\langle
k\rangle^{-m(x)}$ is the ellipticity hypothesis, and the estimates on
$\Delta_k^{\alpha}\partial_x^{\beta}q_0$ follow from those on $p$ by the
quotient rule together with the Leibniz rules for $\Delta_k$ and $\partial_x$.
The class is an algebra in the sense that
$S^{a(\cdot)}_{1,\delta}\cdot S^{b(\cdot)}_{1,\delta} \subseteq
S^{a(\cdot)+b(\cdot)}_{1,\delta}$, again by the Leibniz rules, so an induction on
$N$ using \eqref{eq:parametrix-series} gives $q_N \in
S^{-m(\cdot)-N(1-\delta)}_{1,\delta}$. Because
$\sup_x m(x) < \infty$, the orders $-m(x)-N(1-\delta)$ tend to $-\infty$
uniformly in $x$, so the asymptotic sum may be realised by the usual Borel-type
summation, and the composition formula of \cite{ruzhansky,ruzhansky-book} then
gives $\Op(q)\Op(p)-I$ and $\Op(p)\Op(q)-I$ of order $-\infty$. Every step of
\eqref{eq:parametrix-series} evaluates $p$ and its derivatives at the same base
point $x$, which is the last assertion.
\end{proof}

\subsection{Technical Lemmas}
\label{subsec:technical}

This subsection collects the estimates the main body relies on. The first two
concern the Mat\'ern kernel: its exponential decay, and its expansion near the
origin. We write $\ell>0$ for the range
parameter and $\kappa = \sqrt{2\nu}/\ell$ as in \eqref{eq:lambda}.

\begin{lemma}[Exponential Decay of Mat\'ern Kernel]
\label{lemma:matern-decay}
The Euclidean Mat\'ern kernel with smoothness parameter $\nu > 0$ and range
parameter $\ell > 0$ is
$$k^{\RR}_\nu(r) = \sigma^2 \frac{2^{1-\nu}}{\Gamma(\nu)} \left(\frac{r\sqrt{2\nu}}{\ell}\right)^\nu K_\nu\left(\frac{r\sqrt{2\nu}}{\ell}\right),$$
where $K_\nu$ is the modified Bessel function of the second kind. For every
$j \in \NN_0$ there is $C_j = C_j(\nu,\ell,\sigma)$ with
$$\Bigl|\tfrac{d^j}{dr^j}k^{\RR}_\nu(r)\Bigr| \leq C_j \, r^{\nu-1/2} \, e^{-\kappa r}
\qquad \text{for } r \geq 1 .$$
\end{lemma}

\begin{proof}
Write $z = \kappa r$ with $\kappa = \sqrt{2\nu}/\ell$, so that
$k^{\RR}_\nu(r) = c_{\nu,\sigma}\, z^\nu K_\nu(z)$. The function
$z \mapsto z^\nu K_\nu(z)$ and each of its derivatives satisfy
$\bigl|\tfrac{d^j}{dz^j}\bigl(z^\nu K_\nu(z)\bigr)\bigr| \leq C_j\,
z^{\nu-1/2}e^{-z}$ for $z \geq \kappa$: this is the asymptotic expansion
$K_\mu(z) = \sqrt{\pi/2z}\,e^{-z}\bigl(1+O(1/z)\bigr)$ \cite{dlmf}, valid for
each fixed order $\mu$ and differentiable term by term, together with the
recurrence
$\tfrac{d}{dz}\bigl(z^\mu K_\mu(z)\bigr) = -z^{\mu}K_{\mu-1}(z)$, which shows
that every derivative is again of the same form with $\mu$ decreased by one.
Undoing the substitution multiplies the $j$-th derivative by $\kappa^{j}$ and
absorbs it into $C_j$.
\end{proof}

\begin{lemma}[Local expansion of the Mat\'ern kernel]
\label{lem:matern-expansion}
Let $\nu > 0$. There are entire functions $A$ and $B$ with $B(0) = 1$, and a
constant $c_\nu \neq 0$, such that for $h$ in a neighbourhood of the origin in
$\RR^d$,
\[
    k^{\RR}_\nu(h) \;=\; A\bigl(\|h\|^2\bigr) \;+\;
    c_\nu\,\|h\|^{2\nu}\,B\bigl(\|h\|^2\bigr) \cdot
    \begin{cases}
        1, & \nu \notin \NN, \\
        \log\|h\|, & \nu \in \NN .
    \end{cases}
\]
\end{lemma}

\begin{proof}
Write $z = \kappa\|h\|$, so $k^{\RR}_\nu(h) = c_{\nu,\sigma}\,z^{\nu}K_\nu(z)$,
and recall the ascending series \cite[Eq.~10.25.2]{dlmf}
\[
    I_\mu(z) \;=\; \Bigl(\tfrac z2\Bigr)^{\mu}
        \sum_{j\geq0}\frac{(z^2/4)^{j}}{j!\,\Gamma(\mu+j+1)} .
\]

If $\nu \notin \NN$ then $\nu \notin \ZZ$ and \cite[Eq.~10.27.4]{dlmf} gives
$K_\nu = \tfrac{\pi}{2}\bigl(I_{-\nu}-I_\nu\bigr)/\sin(\nu\pi)$. Multiplying by
$z^{\nu}$,
\[
    z^{\nu}I_{-\nu}(z) = 2^{\nu}\sum_{j\geq0}
        \frac{(z^2/4)^{j}}{j!\,\Gamma(-\nu+j+1)} , \qquad
    z^{\nu}I_{\nu}(z) = 2^{-\nu}z^{2\nu}\sum_{j\geq0}
        \frac{(z^2/4)^{j}}{j!\,\Gamma(\nu+j+1)} ;
\]
the first is an entire function of $z^2$, and the second is $z^{2\nu}$ times an
entire function of $z^2$ whose value at the origin is
$2^{-\nu}/\Gamma(\nu+1) \neq 0$. Substituting $z^2 = \kappa^2\|h\|^2$, collecting the term without the factor
$z^{2\nu}$ into $A$, and absorbing the constants $\kappa^{2\nu}$,
$c_{\nu,\sigma}$, $-\pi/(2\sin\nu\pi)$ and $2^{-\nu}/\Gamma(\nu+1)$ into $c_\nu$
gives the first case; $B$ is then the entire function normalised by $B(0)=1$,
and $c_\nu \neq 0$ because $\sin(\nu\pi)\neq0$.

If $\nu = n \in \NN$, the power series \cite[Eq.~10.31.1]{dlmf} reads
\begin{align*}
    K_n(z) \;=\;& \tfrac12\Bigl(\tfrac z2\Bigr)^{-n}
        \sum_{j=0}^{n-1}\frac{(n-j-1)!}{j!}\Bigl(-\tfrac{z^2}{4}\Bigr)^{j}
      \;-\;(-1)^{n}\log\Bigl(\tfrac z2\Bigr) I_n(z) \\
      &\;+\;(-1)^{n}\tfrac12\Bigl(\tfrac z2\Bigr)^{n}\!\!
        \sum_{j\geq0}\bigl(\psi(j{+}1)+\psi(n{+}j{+}1)\bigr)
        \frac{(z^2/4)^{j}}{j!\,(n+j)!} .
\end{align*}
Multiplying by $z^{n}$: the first sum becomes a polynomial in $z^2$; the third
becomes $z^{2n}$ times an entire function of $z^2$; and
$z^{n}I_n(z) = 2^{-n}z^{2n}\sum_{j\geq0}(z^2/4)^{j}/\bigl(j!\,(n+j)!\bigr)$, so
the middle term is $-(-1)^{n}2^{-n}\bigl(\log z - \log 2\bigr)z^{2n}$ times an
entire function of $z^2$ with value $1/n!$ at the origin. Writing
$\log z = \log\|h\| + \log\kappa$ and collecting all terms without a logarithm
into $A(\|h\|^2)$ gives the second case, with
$c_\nu = -(-1)^{n}2^{-n}\kappa^{2n}c_{\nu,\sigma}/n! \neq 0$.
\end{proof}

The next estimate converts spectral decay into sample-path regularity.

\begin{lemma}[Increment bound for stationary covariances]
\label{lem:increment}
Let $Z$ be a real stationary GeRF on $\TT^d$ whose spectral measure satisfies
$\mu_Z(\{k\}) \leq C_0\,\lambda(k)^{-2s}$ for all $|k| \geq k_0$, where
$s > d/2$. Put $j := \lceil s - d/2\rceil - 1$, so that
$r := 2(s-d/2) - 2j \in (0,2]$. Then $\rho_Z \in C^{2j}(\TT^d)$, and for every
multi-index $\alpha$ with $|\alpha| = 2j$ and every $\epsilon < r/2$,
\[
    \bigl|\partial^\alpha \rho_Z(h) - \partial^\alpha \rho_Z(0)\bigr|
    \;=\; \calO\bigl(\|h\|^{2\epsilon}\bigr), \qquad h\to0 .
\]
Consequently, by the ``if'' half of Theorem~\ref{thm:criterion}, $Z$ has
samples in $C^{(s-d/2)^-}_{loc}(\calO)$.
\end{lemma}

\begin{proof}
Only the tail matters. The finitely many frequencies $|k| < k_0$ contribute to
$\rho_Z$ a trigonometric polynomial, which is smooth and even, so each of its
even-order derivatives has vanishing gradient at the origin and its increment
there is $\calO(\|h\|^2)$; since $2\epsilon < r \leq 2$ this is absorbed into the
claim. We may therefore assume $\mu_Z(\{k\}) \leq C_0\lambda(k)^{-2s}$ for all
$k$.

Write $\rho_Z(h) = \sum_{k} \mu_Z(\{k\})\, e^{2\pi i k\cdot h}$. Because
$|\alpha| = 2j$ and $2j - 2s + d < 0$, the differentiated series
\begin{equation}\label{eq:diff-series}
    \partial^\alpha \rho_Z(h)
    \;=\; \sum_{k} (2\pi i k)^{\alpha}\, \mu_Z(\{k\})\, e^{2\pi i k\cdot h}
\end{equation}
converges absolutely and uniformly, whence $\rho_Z \in C^{2j}$. Since $Z$ is
real and stationary, $\mu_Z$ is even; pairing $k$ with $-k$ in
\eqref{eq:diff-series} and using that $|\alpha|$ is even replaces the
exponential by a cosine, so
\[
    \partial^\alpha \rho_Z(h) - \partial^\alpha \rho_Z(0)
    \;=\; \sum_{k} (2\pi i k)^{\alpha}\, \mu_Z(\{k\})
          \bigl(\cos(2\pi k\cdot h) - 1\bigr) .
\]
For $\theta \in [0,2]$, interpolating $|\cos t - 1| \leq 2$ against
$|\cos t - 1| \leq \tfrac12 t^2$ gives
$|\cos t - 1| \leq 2^{1-\theta}|t|^{\theta}$, so
$|\cos(2\pi k\cdot h) - 1| \leq 2^{1-\theta}(2\pi)^{\theta}|k|^{\theta}\|h\|^{\theta}$
and
\[
    \bigl|\partial^\alpha \rho_Z(h) - \partial^\alpha \rho_Z(0)\bigr|
    \;\leq\; C_\theta\, \|h\|^{\theta} \sum_{k} \langle k\rangle^{\,2j - 2s + \theta} .
\]
The series converges precisely when $2j - 2s + \theta + d < 0$, that is when
$\theta < 2(s - d/2) - 2j = r$. Taking $\theta = 2\epsilon$ gives the display,
and $j + \epsilon$ ranges up to $j + r/2 = s - d/2$.
\end{proof}

The last estimate is the two-sided lattice bound used in the proof of
Theorem~\ref{thm:pointwise}. It is stated with a weight because the induction
on the order of differentiation there produces the multipliers
\eqref{eq:theta}.

\begin{lemma}[Weighted lattice estimate]\label{lem:lattice}
Let $w : \ZZ^d \to [0,\infty)$ satisfy $w \leq C_1$ everywhere and
$w(k) \geq c_1 > 0$ for $|k| \geq k_1$, and write
\[
    S_s(h) \;:=\; \sum_{k\in\ZZ^d} w(k)\,\lambda(k)^{-2s}
        \bigl(1-\cos 2\pi k\cdot h\bigr) .
\]
\begin{enumerate}
\item[(a)] If $d < 2s < d+2$ then there are constants $0<c\leq C<\infty$ with
\[
    c\,|h|^{2s-d} \;\leq\; S_s(h) \;\leq\; C\,|h|^{2s-d}
    \qquad \text{for } 0 < |h| \leq \tfrac14 .
\]
The constants depend only on $d,c_1,C_1,k_1,\kappa$ and on $s$, and may be
chosen uniformly for $s$ ranging in a compact subinterval of
$\bigl(\tfrac d2,\tfrac{d+2}2\bigr)$.
\item[(b)] If $2s = d+2$ then $S_s(h) \leq C\,|h|^{2}\log(2/|h|)$ for
$0<|h|\leq\tfrac14$. In particular $S_s(h) = \calO(|h|^{2\epsilon})$ for every
$\epsilon<1$.
\end{enumerate}
For $w \equiv 1$, part (a) is the periodic form of the classical Bessel
potential asymptotics; see \cite{stein}.
\end{lemma}

\begin{proof}
Throughout $\lambda(k) \asymp \langle k\rangle$ with constants independent of
$s$ once $s$ is confined to a bounded set, so it suffices to prove the bounds
with $\langle k\rangle^{-2s}$ in place of $\lambda(k)^{-2s}$. Write $S(h)$ for
$S_s(h)$.

\emph{Upper bound, case \textup{(a)}.} Using $w \leq C_1$, split at $|k| \leq 1/|h|$. For the near
frequencies, $1-\cos t \leq t^2/2$ gives
\[
    \sum_{|k|\leq 1/|h|}\!\!\langle k\rangle^{-2s}\bigl(1-\cos2\pi k\cdot h\bigr)
    \;\leq\; 2\pi^2|h|^2\!\!\sum_{|k|\leq 1/|h|}\!\!\langle k\rangle^{2-2s}
    \;\leq\; C|h|^2\,|h|^{\,2s-d-2} \;=\; C|h|^{2s-d},
\]
where $\sum_{|k|\leq R}\langle k\rangle^{2-2s} \leq CR^{\,d+2-2s}$ because
$d+2-2s>0$. For the far frequencies, $1-\cos \leq 2$ gives
$\sum_{|k|> 1/|h|}\langle k\rangle^{-2s} \leq C|h|^{2s-d}$ because $2s>d$. Both
constants stay bounded as long as $d+2-2s$ and $2s-d$ stay bounded away from
$0$, which is the asserted uniformity in $s$.

\emph{Upper bound, case \textup{(b)}.} If $2s = d+2$ the far sum is unchanged
and contributes $\calO(|h|^{2})$, while the near sum becomes
$|h|^2\sum_{|k|\leq1/|h|}\langle k\rangle^{-d} \leq C|h|^2\log(2/|h|)$; adding
the two gives the claim, and $\log(2/|h|) = \calO(|h|^{-2(1-\epsilon)})$ for
every $\epsilon<1$.

\emph{Lower bound, case \textup{(a)}.} Fix $h_0 := \min\bigl(\tfrac14,\tfrac{1}{8k_1},h_1\bigr)$,
where $h_1$ is chosen below, and assume first $0<|h| \leq h_0$. Let
\[
    A_h \;:=\; \Bigl\{k \in \ZZ^d \;:\; \tfrac{1}{8|h|} \leq |k| \leq \tfrac{1}{4|h|}\Bigr\} .
\]
Choosing $h_1$ small enough that the inner radius exceeds an absolute constant,
the lattice-point count in this annulus satisfies $\#A_h \geq c_2|h|^{-d}$, and
every $k \in A_h$ has $\langle k\rangle|h| \in [c_4,C_4]$ with $c_4,C_4$
absolute, so that $\langle k\rangle^{-2s} \geq c_3|h|^{2s}$ with
$c_3 = \min(c_4,C_4)^{-2s}\wedge\max(c_4,C_4)^{-2s}$ bounded below uniformly for
bounded $s$; and, since
$|k| \geq 1/(8|h|) \geq k_1$, also $w(k) \geq c_1$. For $k \in A_h$ we have
$|2\pi k\cdot h| \leq 2\pi|k||h| \leq \pi/2$, so the elementary inequality
$1-\cos t \geq \tfrac{2}{\pi^2}t^2$ for $|t| \leq \pi$ gives
$1-\cos 2\pi k\cdot h \geq 8\,(k\cdot h)^2$. Hence
\[
    S(h) \;\geq\; 8\,c_1c_3\,|h|^{2s} \sum_{k \in A_h}(k\cdot h)^2
    \;=\; 8\,c_1c_3\,|h|^{2s}\; h^{\top}\!M_h\,h ,
    \qquad M_h := \sum_{k\in A_h}k\,k^{\top} .
\]
The set $A_h$ is defined by $|k|$ alone, so it is invariant under every signed
permutation of the coordinates. Writing $P$ for the orthogonal matrix of such a
transformation, $PA_h = A_h$ and therefore $PM_hP^{\top} = M_h$. A matrix
commuting with all sign changes has vanishing off-diagonal entries, and one
commuting with all coordinate permutations has equal diagonal entries; hence
$M_h = \tfrac1d\operatorname{tr}(M_h)\,I = \tfrac1d\bigl(\sum_{k\in A_h}|k|^2\bigr)I$
exactly. Consequently
\[
    h^{\top}M_h\,h \;=\; \frac{|h|^2}{d}\sum_{k\in A_h}|k|^2
    \;\geq\; \frac{|h|^2}{d}\cdot\frac{1}{64|h|^2}\cdot c_2|h|^{-d}
    \;=\; \frac{c_2}{64d}\,|h|^{-d} ,
\]
and $S(h) \geq c\,|h|^{2s-d}$ as claimed.

For $h_0 \leq |h| \leq \tfrac14$ the sum $S$ converges uniformly, since
$S(h) \leq 2C_1\sum_k\langle k\rangle^{-2s} < \infty$, and is therefore
continuous; it is also strictly positive there, for $S(h)=0$ would force
$k\cdot h \in \ZZ$ for every $k$ with $|k| \geq k_1$ and hence $h \in \ZZ^d$,
which is impossible for $0<|h|\leq\frac14$ on $\TT^d$. A continuous positive
function on the compact set $\{h_0 \leq |h| \leq \frac14\}$ is bounded below,
and $|h|^{2s-d}$ is bounded above there, so the stated inequality persists after
enlarging the constant. Both bounds are joint continuous functions of $(s,h)$ on
the compact set obtained by letting $s$ range over a compact subinterval of
$(\tfrac d2,\tfrac{d+2}2)$, which gives the uniformity in $s$.
\end{proof}

\subsection{Periodisation and truncation}
\label{subsec:periodisation}

The covariances used in this paper are periodisations of Euclidean ones, and two
of the arguments above rely on the fact that periodisation is invisible near the
diagonal.

Throughout, $k^{\RR}_\nu$ is the Euclidean Mat\'ern kernel of
Lemma~\ref{lemma:matern-decay} and
\begin{equation}\label{eq:periodisation}
    k_\nu(h) \;=\; \sum_{n \in \ZZ^d} k^{\RR}_\nu\bigl(\|h+n\|\bigr),
    \qquad h \in \TT^d ,
\end{equation}
is its periodisation, the covariance \eqref{matern-cov}.

\begin{proposition}[Periodisation is smooth near the diagonal]
\label{prop:periodisation}
Let $R_\nu(h) := k_\nu(h) - k^{\RR}_\nu(\|h\|)$ denote the sum of the
off-lattice terms in \eqref{eq:periodisation}. Then $R_\nu \in C^\infty$ on
$\{\|h\|_\infty < \tfrac14\}$, and for every multi-index $\alpha$ there is
$C_\alpha = C_\alpha(\nu,\ell,\sigma,d)$ with
\[
    \sup_{\|h\|_\infty \leq \tfrac14}\bigl|\partial^\alpha R_\nu(h)\bigr|
    \;\leq\; C_\alpha \, e^{-\kappa/(2\sqrt d)} .
\]
The same bounds hold on a complex neighbourhood of $\{\|h\|_\infty<\tfrac14\}$,
so $R_\nu$ is real-analytic there. These are the bounds the increment comparison
in the proof of Proposition~\ref{proposition:regularity} uses to carry
Proposition~\ref{euclidean-regularity} over to the torus.
\end{proposition}

\begin{proof}
Let $\|h\|_\infty \leq \tfrac14$ and $n \in \ZZ^d\setminus\{0\}$, so
$\|n\|_\infty \geq 1$. Then
$\|h+n\|_\infty \geq \|n\|_\infty - \tfrac14 \geq \tfrac34\|n\|_\infty$, whence
\begin{equation}\label{eq:lattice-sep}
    \|h+n\| \;\geq\; \|h+n\|_\infty \;\geq\; \tfrac34\|n\|_\infty
    \;\geq\; \tfrac{3}{4\sqrt d}\,\|n\| \;\geq\; \tfrac{3}{4\sqrt d} .
\end{equation}
In particular no term with $n \neq 0$ has its singularity in the region
considered, and $h \mapsto k^{\RR}_\nu(\|h+n\|)$ is smooth there. Differentiating
$|\alpha|$ times produces derivatives of $k^{\RR}_\nu$ of order at most
$|\alpha|$ multiplied by derivatives of $h \mapsto \|h+n\|$, which are bounded
uniformly by \eqref{eq:lattice-sep}; Lemma~\ref{lemma:matern-decay} and
\eqref{eq:lattice-sep} then give
\[
    \bigl|\partial^\alpha_h\, k^{\RR}_\nu(\|h+n\|)\bigr|
    \;\leq\; C_\alpha' \, \|n\|^{\nu-1/2}\,
             e^{-3\kappa\|n\|/(4\sqrt d)} .
\]
Splitting the exponential as
$e^{-3\kappa\|n\|/(4\sqrt d)} = e^{-\kappa\|n\|/(2\sqrt d)}\,
e^{-\kappa\|n\|/(4\sqrt d)} \leq e^{-\kappa/(2\sqrt d)}\,e^{-\kappa\|n\|/(4\sqrt d)}$,
which uses $\|n\|\geq1$, and summing the convergent series
$\sum_{n\neq0}\|n\|^{\nu-1/2}e^{-\kappa\|n\|/(4\sqrt d)}$ gives the stated bound
after relabelling constants. Uniform convergence of the
differentiated series on the region gives $R_\nu \in C^\infty$. For
analyticity, $z \mapsto z^\nu K_\nu(z)$ is holomorphic on $\Re z > 0$ and the
same estimates hold with $\|h+n\|$ replaced by its holomorphic extension on a
complex neighbourhood on which \eqref{eq:lattice-sep} persists, so the series
converges uniformly there and $R_\nu$ is the restriction of a holomorphic
function.
\end{proof}

\begin{corollary}[Truncation error]
\label{cor:truncation-error}
Truncating \eqref{eq:periodisation} at $\|n\|_\infty \leq N$ leaves an error
\[
    \sup_{h\in\TT^d}\Bigl|k_\nu(h) - \sum_{\|n\|_\infty \leq N} k^{\RR}_\nu(\|h+n\|)\Bigr|
    \;\leq\; C\,N^{\,d-1+\nu-1/2}\, e^{-3\kappa N/(4\sqrt d)}
\]
for $N \geq 1$, with $C$ depending on $\nu$, $\ell$, $\sigma$ and $d$ only. The
same bound with $C_\alpha$ in place of $C$ holds for each derivative
$\partial^\alpha k_\nu$.
\end{corollary}

\begin{proof}
By \eqref{eq:lattice-sep} and Lemma~\ref{lemma:matern-decay} each term with
$\|n\|_\infty = M$ is at most
\[
    C\,M^{\nu-1/2}\,e^{-3\kappa M/(4\sqrt d)} ,
\]
and there are $O(M^{d-1})$ such lattice points. Summing the resulting series over
$M > N$ and bounding it by its first term times a convergent geometric factor
gives the claim.
\end{proof}

\subsection{Proof of Proposition~\ref{proposition:regularity}}
\label{app:proof-regularity}

\begin{proof}[Proof of Proposition~\ref{proposition:regularity}]
Write $\rho_Z(x-x') = k_\nu(x,x')$, so that the stationary covariance
\eqref{matern-cov} defines the operator $L_{\rho_Z}f = \rho_Z \ast f$ on
$C^\infty(\TT^d)$. Set
\[
    j \;:=\; \lceil\nu\rceil-1 , \qquad r \;:=\; 2\nu-2j \in (0,2] ,
\]
so that $\nu = j + r/2$, with $r<2$ exactly when $\nu \notin \NN$.

By Poisson summation the Fourier coefficients of a periodisation are the samples
at integer frequencies of the Euclidean Fourier transform of the summand: the
identity $\sum_{n \in \ZZ^d} f(h+n) = \sum_{k \in \ZZ^d} \hat f(k)e^{2\pi i
k\cdot h}$ holds pointwise, with both series absolutely convergent, as soon as
$|f(x)| \leq C(1+|x|)^{-d-\epsilon}$ and $|\hat f(\xi)| \leq
C(1+|\xi|)^{-d-\epsilon}$ for some $\epsilon>0$. Both hypotheses hold for the
Euclidean Mat\'ern kernel $k^{\RR}_\nu$ of Lemma~\ref{lemma:matern-decay}: that
lemma gives exponential decay of $k^{\RR}_\nu$ itself, and its transform is the
Whittle--Mat\'ern spectral density of \cite{lindgren2011},
\begin{equation}\label{eq:whittle-matern}
    \widehat{k^{\RR}_\nu}(\xi)
    \;=\; c_1\Bigl(\tfrac{2\nu}{\ell^2} + 4\pi^2|\xi|^2\Bigr)^{-(\nu + d/2)},
    \qquad c_1 = c_1(\nu,\ell,\sigma,d) > 0 ,
\end{equation}
which is $\calO(|\xi|^{-2\nu-d})$ and therefore integrable precisely because
$\nu>0$.

Convolution on $\TT^d$ is a Fourier multiplier, $\widehat{\rho \ast f}(k) =
\hat\rho(k)\hat f(k)$, and
the previous paragraph identifies $\hat\rho_Z(k) = \widehat{k^{\RR}_\nu}(k)$.
Since $\kappa^2 = 2\nu/\ell^2$, the operator $L_{\rho_Z}$ is therefore a Fourier
multiplier with symbol
\begin{equation}\label{eq:correct-symbol}
    \sigma(k) \;=\; c_1\bigl(\kappa^2 + 4\pi^2|k|^2\bigr)^{-(\nu+d/2)}
    \;=\; c_1\,\lambda(k)^{-2s}, \qquad s := \nu + \tfrac d2 ,
\end{equation}
of order $m = -2\nu-d$.

Lemma~\ref{lem:increment} applied with this $s$ gives
\[
    \bigl|\partial^\alpha \rho_Z(h) - \partial^\alpha \rho_Z(0)\bigr| = \calO(\|h\|^{2\epsilon}),
    \qquad |\alpha| = 2j, \quad \epsilon < r/2 ,
\]
so the stationary criterion \eqref{eq:criterion-stat} holds with $n=j$ and
$\gamma = r/2$. By the ``if'' half of Theorem~\ref{thm:criterion}, samples
lie in $C^{(j+r/2)^-}_{loc} = C^{\nu^-}_{loc}$.

Proposition~\ref{prop:periodisation} writes $k_\nu = k^{\RR}_\nu + R_\nu$ with
$R_\nu$ real-analytic on $\{\|h\|_\infty<\tfrac14\}$, and
Lemma~\ref{lem:matern-expansion} expands the Euclidean kernel near the origin as
\begin{equation}\label{eq:matern-local}
    k^{\RR}_\nu(h) \;=\; A\bigl(\|h\|^2\bigr) \;+\;
    \begin{cases}
        c_\nu\,\|h\|^{2\nu}\,B\bigl(\|h\|^2\bigr), & \nu \notin \NN, \\[3pt]
        c_\nu\,\|h\|^{2\nu}\log\|h\|\;B\bigl(\|h\|^2\bigr), & \nu \in \NN,
    \end{cases}
\end{equation}
with $A$ and $B$ entire, $B(0)=1$ and $c_\nu \neq 0$; the logarithm in the
second line is the one carried by $K_\nu$ at integer order.

Fix $\alpha := 2j\,e_1$, so $|\alpha| = 2j$, and take $h = t\,e_1$ with
$t \downarrow 0$. The contributions of $A(\|h\|^2)$ and of $R_\nu$ to
$\partial^{\alpha}\rho_Z(h) - \partial^{\alpha}\rho_Z(0)$ are $\calO(t^2)$: both
functions are even and $C^\infty$ near the origin, so $\partial^\alpha$ of each
is again even and $C^\infty$ there, hence has vanishing gradient at $0$, and
Taylor's theorem applies on the segment from $0$ to $h$, which lies in
$\{\|h\|_\infty<\tfrac14\}$. By
Proposition~\ref{prop:periodisation},
\begin{equation}\label{eq:remainder-increment}
    \bigl|\partial^\alpha R_\nu(h) - \partial^\alpha R_\nu(0)\bigr|
    \;\leq\; C'_\alpha\, e^{-\kappa/(2\sqrt d)}\,\|h\|^2 ,
    \qquad \|h\|_\infty \leq \tfrac14 ,
\end{equation}
with $C'_\alpha$ depending on $\alpha$ and on $\nu,\ell,\sigma,d$ (Taylor's
remainder uses the bounds of Proposition~\ref{prop:periodisation} at order
$|\alpha|+2$).
For the singular term, $\partial_1^{2j}|t|^{2\nu} = \gamma_\nu\,|t|^{\,r}$ with
\[
    \gamma_\nu \;=\; 2\nu(2\nu-1)\cdots(2\nu-2j+1) \;\neq\; 0 ,
\]
the factors being strictly positive because $2\nu > 2j$; the corrections
supplied by $B$ are $\calO(|t|^{\,r+2})$. Hence
\begin{equation}\label{eq:exact-increment}
    \partial^{\alpha}\rho_Z(t e_1) - \partial^{\alpha}\rho_Z(0)
    \;=\;
    \begin{cases}
        c_\nu\gamma_\nu\,|t|^{\,r}\bigl(1+\calO(t^2)\bigr) + \calO(t^2),
            & \nu \notin \NN, \\[3pt]
        c_\nu\gamma'_\nu\,t^{2}\log(1/|t|)\bigl(1+o(1)\bigr),
            & \nu \in \NN,
    \end{cases}
\end{equation}
with $\gamma'_\nu \neq 0$. When $\nu \notin \NN$ we have $r < 2$, so the first
line is of exact order $|t|^{\,r}$; when $\nu \in \NN$ the logarithm dominates
every $\calO(t^2)$ contribution.

Let $\nu' > \nu$ and suppose samples lay in
$C^{\nu'^-}_{loc}$; write $\nu' = n' + \gamma'$ with $\gamma' \in (0,1]$.

If $\nu \notin \NN$, then $\nu = j + r/2$ with $r/2 \in (0,1)$, and since
$C^{\nu'^-}_{loc} \subseteq C^{\nu''^-}_{loc}$ for $\nu \le \nu'' \le \nu'$ we
may take $\nu'' = j+\gamma''$ with $\gamma'' \in \bigl(r/2,\min(1,\nu'-j)\bigr)$.
Choosing $\epsilon \in (r/2,\gamma'')$, the ``only if'' half of
Theorem~\ref{thm:criterion} would force
$|\partial^\alpha\rho_Z(h)-\partial^\alpha\rho_Z(0)| = \calO(\|h\|^{2\epsilon})$
with $2\epsilon > r$, contradicting the first line of
\eqref{eq:exact-increment}.

If $\nu \in \NN$, then $r = 2$ and $j = \nu-1$, so $\gamma' \leq 1$ forces
$n' \geq \nu$. The criterion at order $n'$ requires $\rho_Z \in C^{2n'}
\subseteq C^{2\nu}$. But differentiating \eqref{eq:matern-local} $2\nu$ times in
the direction $e_1$ produces a term $c_\nu\gamma''_\nu\log|t|$ with
$\gamma''_\nu \neq 0$, which is unbounded as $t\to0$, while every other
contribution stays bounded; hence $\rho_Z \notin C^{2\nu}(\TT^d)$ and the
criterion is not applicable at order $n'$.

Either way the assumption fails, and the samples lie in $C^{\nu^-}_{loc}$ and no
more.
\end{proof}

\end{appendices}

\end{document}